\documentclass[12pt,a4paper]{article}

\usepackage{amsmath,amsthm,amssymb}
\usepackage{hyperref}
\usepackage{tikz-cd}
\usepackage{notespackage}
\usepackage{cleveref}
\usepackage{geometry}
\title{Tian--Yau-type metrics from K3 divisors with ADE contractions}
\author{Yujian Song}
\date{\today}

\begin{document}

\maketitle

\begin{abstract}
	We construct complete Ricci-flat K\"{a}hler metrics on the complement of a smooth anti-canonical K3 divisor in a compact K\"{a}hler threefold, assuming that the normal bundle of the divisor is nef and big. Such a metric exists in every K\"{a}hler class whose exceptional periods are positive. We use a generalized Gibbons--Hawking ansatz to construct an approximate model on the normal bundle. When the normal bundle contracts an ADE configuration, this model combines an expanding Ricci-flat orbifold K3 metric with ALE metrics whose exceptional periods remain fixed. After improving the model by a finite iteration and transferring it to the divisor complement, we apply the Tian--Yau--Hein existence theorem.
\end{abstract}

\tableofcontents

\section{Introduction}\label{intro}

The study of Ricci-flat K\"{a}hler metrics is a central problem in K\"{a}hler geometry. Yau's solution of the Calabi conjecture \cite{yauRicciCurvatureCompact1978} establishes the existence of Ricci-flat K\"{a}hler metrics in every K\"{a}hler class on a compact K\"{a}hler manifold with vanishing first real Chern class. In the noncompact setting, the existence problem also depends on the geometry at infinity. The foundational work of Tian and Yau \cite{tianCompleteKahlerManifolds1990,tianCompleteKahlerManifolds1991} produced complete Ricci-flat K\"{a}hler metrics on broad classes of quasi-projective manifolds, including complements of smooth anti-canonical divisors in Fano manifolds. Their approach, subsequently developed by Hein \cite{heinGravitationalInstantons2010,heinGravitationalInstantonsRational2012}, provides a framework for solving the complex Monge--Amp\`{e}re equation on a noncompact manifold once a suitable complete background metric has been constructed. We refer to this analytic framework as the Tian--Yau--Hein package.

For the classical Tian--Yau construction with a smooth anti-canonical divisor, the geometry at infinity is described by the Calabi model on the punctured normal bundle of the divisor. The ampleness of the normal bundle supplies a Ricci-flat K\"{a}hler metric on the divisor in the corresponding positive class, from which the model is obtained by the Calabi ansatz. The geometry and analysis of this model play an important role in the work of Hein--Sun--Viaclovsky--Zhang on collapsing K3 metrics \cite{heinNilpotentStructuresCollapsing2021}. Their subsequent work \cite{heinAsymptoticallyCalabiMetrics2025} relates asymptotically Calabi metrics to weak Fano compactifications and a generalized Tian--Yau construction. More recently, Chen \cite{chenCalabiYauMetricsCalabi2024} constructed Calabi-type metrics in prescribed cohomology classes, with polynomial convergence to the Calabi model, assuming that the smooth anti-canonical divisor has ample normal bundle. Thus, even when the ambient Fano hypothesis is relaxed, positivity of the normal bundle remains an important ingredient in these constructions.

Related existence results exhibit other geometries at infinity. Haskins--Hein--Nordstr\"{o}m \cite{haskinsAsymptoticallyCylindricalCalabi2015} developed the theory of asymptotically cylindrical Calabi--Yau manifolds, while Conlon--Hein \cite{conlonAsymptoticallyConicalCalabi2015} constructed asymptotically conical metrics in related divisor complement settings. Collins--Li \cite{collinsCompleteCalabiYau2024} treated complements of two transversely intersecting divisors forming an anti-canonical divisor, where a different model is required. Further examples include Li's complete Calabi--Yau metric on $\mathbb C^3$ \cite{liNewCompleteCalabi2019}, Wang's metrics of generalized ALG asymptotics \cite{wangRicciflatManifoldsGeneralized2023}, and Min's higher-dimensional ALF constructions \cite{minConstructionHigherDimensional2023}. These results illustrate the importance of finding a model adapted to the particular geometry of the compactification.

In this paper, we consider smooth anti-canonical K3 divisors in compact K\"{a}hler threefolds and replace ampleness of the normal bundle by the weaker assumption that it is nef and big. When the normal bundle is not ample, its first Chern class has zero intersection with certain curves on the divisor and therefore cannot be represented by a smooth K\"{a}hler form. The usual Calabi model is consequently unavailable. The K3 setting provides a replacement: the curves of zero intersection can be contracted to ADE singularities, and the resulting orbifold geometry can be combined with Ricci-flat ALE metrics near the exceptional locus. The construction must retain the exceptional periods prescribed by the K\"{a}hler class on the complement.

More precisely, let $\overline{M}$ be a smooth compact K\"{a}hler threefold, $D\in |-K_{\overline{M}}|$ be a smooth anti-canonical K3 divisor, $M=\overline{M}\setminus D$. The normal bundle $L=\mathcal{O}_{\overline{M}}(D)|_D$ is nef and big. Let $\Omega_{M}$ be a meromorphic volume form with $\operatorname{div}\Omega_{M}=-D$.

For sufficiently large $k\gg1$, $L^{\otimes k}$ is base-point-free, and induces a morphism \begin{equation}
	\pi_Y\colon D\to Y\subset \mathbb{P}H^0(D,L^{\otimes k})
\end{equation} that contracts precisely the curves orthogonal to $L$. The surface $Y$ has only ADE singularities, and $D$ is its minimal crepant resolution; see \cite[Chapter2, Remark 3.4]{huybrechtsLecturesK3Surfaces2016}. Denote the irreducible exceptional components by $E_1,\ldots,E_N$.

To formulate the condition on the prescribed class, let $j:M\hookrightarrow\overline{M}$ be the inclusion. By the Thom--Gysin sequence \begin{equation}
	\begin{tikzcd}
		H^0(D,\mathbb{R}) \ar[r,"{1\mapsto [D]}"] & H^2(\overline{M},\mathbb{R}) \ar[r,"j"] & H^2(M,\mathbb{R}) \ar[r] & H^1(D,\mathbb{R})=0,
	\end{tikzcd}
\end{equation} every cohomology class $\mathfrak{c}\in H^2(M,\mathbb{R})$ has a lifting $\bar{\mathfrak{c}}\in H^2(\overline{M},\mathbb{R})$, and two liftings differ by a multiple of $[D]$. Since $[D]\cdot E_i=0$, the intersection $p_i(\mathfrak{c})=\bar{\mathfrak{c}}\cdot E_i$ is independent of the lifting.

\begin{definition}
	We call $p_i(\mathfrak{c})=\bar{\mathfrak{c}}\cdot E_i$ the exceptional period of $\mathfrak{c}$ at $E_i$.
\end{definition}

Our main theorem is the following.

\begin{theorem}
	Let $\overline{M}$, $D$, $M$, $L$, and $\Omega_M$ be as above. For every K\"{a}hler class $\mathfrak{c}\in H^2(M,\mathbb{R})$ satisfying \begin{equation}
		p_i(\mathfrak{c})>0,\quad i=1,\ldots,N,
	\end{equation} there exists a complete Ricci-flat K\"{a}hler metric $\omega_{\mathrm{CY}}$, such that \begin{equation}
		[\omega_{\mathrm{CY}}]=\mathfrak{c},\quad\omega_{\mathrm{CY}}^3=6\sqrt{-1}\Omega_{M}\wedge\overline{\Omega}_{M}.
	\end{equation}
\end{theorem}

When $L$ is ample, there are no exceptional components and the period condition is vacuous. The theorem therefore extends the prescribed-class existence statement to the nef and big setting in complex dimension three, subject to positivity of the exceptional periods. These periods determine the areas of the exceptional curves in the ALE pieces used in the model construction.

The geometric ingredients are closely related to the study of degenerating Ricci-flat metrics on K3 surfaces. Tosatti \cite{tosattiLimitsCalabiYauMetrics2009} proved smooth convergence away from the exceptional set for Ricci-flat metrics whose K\"ahler classes approach a nef and big boundary class on a projective Calabi--Yau manifold. In complex dimension two, orbifold degenerations and their resolutions are reflected in Kobayashi's work \cite{kobayashiModuliEinsteinMetrics} and the weighted analytic treatments of Jiang \cite{jiang2026kummerconstructioncalabiyauhyperkahler}. The local building blocks for general ADE singularities are provided by Kronheimer's ALE hyper-K\"ahler spaces \cite{kronheimerConstructionALESpaces1989}; Joyce's analysis \cite{joyceAsymptoticallyLocallyEuclidean2001} supplies the prescribed-class Ricci-flat ALE metrics and weighted elliptic theory used here. On the compact orbifold side, we use the orbifold Calabi--Yau theorem; see \cite{faulkYausTheoremEffective2019}.

Our task is not simply to desingularize a compact orbifold metric. We must construct a family of metrics on the divisor whose dependence on a radial parameter is compatible with a K\"ahler metric on its punctured normal bundle. For this purpose, we use the generalized Gibbons--Hawking ansatz in the form described by Sun--Zhang \cite[Section~2]{sunComplexStructureDegenerations2020}. The circle action on the fibers reduces the problem to a family of K\"ahler forms $\omega_t$ on $D$, together with a positive function $h_t$, where $t$ is the moment map coordinate. After choosing a suitable closed real $(1,1)$-form $\beta$ on $\overline M$ with $[\beta|_M]=c$, the reduced classes satisfy \begin{equation}
	[\omega_t]=2\pi tc_1(L)+[\beta|_D].
\end{equation} In particular, \begin{equation}
	\int_{E_i}\omega_t=p_i(\mathfrak{c})
\end{equation} is independent of $t$.

We construct the family $\omega_t$ by gluing an expanding Ricci-flat orbifold K3 metric to ALE metrics with these fixed exceptional periods. Away from the exceptional locus, its leading behavior is $t\pi_Y^*\omega_Y$, where $\omega_Y$ is the orbifold Ricci-flat metric in the class descended from $2\pi c_1(L)$. Near the exceptional curves, the model retains the ALE geometry at the scale determined by the prescribed periods. After rescaling by $t^{-1}$, the exceptional areas shrink to zero, whereas before rescaling they remain fixed. This coexistence of expanding orbifold geometry and fixed exceptional areas is the feature that replaces the homothetic family in the classical Calabi ansatz.

The initial gluing does not give sufficiently small error for a direct application of the Tian--Yau--Hein package. Moreover, estimates on individual metrics in the family are not enough, since the reduced equations also involve derivatives with respect to $t$. We introduce weighted H\"{o}lder spaces adapted to the orbifold and ALE regions, construct a uniformly bounded right inverse for the Laplacian, and establish estimates for its parameter dependence. A finite iteration then improves the approximate solution while preserving the reduced cohomology classes and the positivity required by the ansatz.

Finally, we reconstruct a K\"{a}hler metric on a punctured neighborhood of the zero section in $L$ and transfer it to the end of $M$. After gluing it to a global representative of the prescribed class, we adjust a parameter in the model to satisfy the integral compatibility condition for the complex Monge--Amp\`{e}re equation. The resulting complete background metric has the volume growth, local geometric control, and Ricci potential decay required by the Tian--Yau--Hein package. Solving the Monge--Amp\`{e}re equation produces a complete Ricci-flat K\"{a}hler metric in the original class $c$.

The paper is organized as follows. \Cref{s2} develops the generalized Gibbons--Hawking ansatz and reconstructs the global potential from the reduced data. \Cref{s3} treats the cohomological preparation, constructs the approximate family by orbifold-ALE gluing, and introduces the weighted spaces. \Cref{s4} establishes the uniform inverse estimates for the Laplacian. \Cref{s5} studies parameter dependence and carries out the finite improvement procedure. \Cref{s6} transfers the model to the divisor complement and imposes the integral condition. \Cref{s7} verifies the remaining geometric hypotheses and completes the proof.

Throughout the paper, we use \begin{equation}
	d^c=Jd=\sqrt{-1}(\bar{\partial}-\partial),\quad dd^c=2\sqrt{-1}\partial\bar{\partial},
\end{equation} and the Laplacian convention \begin{equation}
	\Delta_\omega u=\frac{n\sqrt{-1}\partial\bar{\partial}u\wedge\omega^{n-1}}{\omega^n}.
\end{equation}

\noindent\textbf{Acknowledgment}: The author is grateful to Prof. Gang Tian for his guidance and encouragement. The author also thanks Ruiming Liang for helpful discussions.

\section{Generalized Gibbons--Hawking ansatz}\label{s2}

This section establishes the reduction and reconstruction framework used to
build the metric on the normal bundle. The generalized Gibbons--Hawking
ansatz expresses an $S^1$-invariant metric in terms of a family of K\"{a}hler
forms on the base and a positive function, parametrized by the moment map
coordinate $t$. We distinguish the compatibility conditions needed to
reconstruct a K\"{a}hler metric from the additional volume form identity
imposing the Calabi--Yau condition. This distinction allows us to work with
approximate solutions in the subsequent sections. We also recover a global
potential relative to a fixed closed $(1,1)$-form, which will be used to
transfer the model to the divisor complement in \cref{s6}.

\subsection{K\"{a}hler reduction of $S^1$-symmetry}

Let $X^n$ be a Calabi--Yau manifold of dimension $n$, $\pi_L:L\to X$ be a nef and big line bundle. Let $N$ be a small punctured tubular neighborhood of the zero section in $L$. We introduce the generalized Gibbons--Hawking ansatz on $N$. The following discussion is essentially in \cite[section 2]{sunComplexStructureDegenerations2020}, we present it here for completeness.

There is a natural $S^1$-action on $L^\times$ by rotating the fibers. Assume there is a K\"{a}hler metric $\omega_L$ on $N$ such that the $S^1$-action is holomorphic, isometric, and Hamiltonian. Let $\xi$ be a Killing field generating the $S^1$-action. Since the action is Hamiltonian, there is a moment map $t$ such that \begin{equation}\label{1}
	dt=\iota_\xi\omega_L.
\end{equation}

Denote by $\{z_1,\ldots,z_n\}$ the local holomorphic coordinates on $X$. Then they can be viewed as local holomorphic functions on $N$. Let the fiber coordinate $w=|w|e^{\sqrt{-1}\sigma}$. Then $\{t,\sigma,z_1,\ldots,z_n\}$ gives a local coordinate system on $N$, and we have $\xi=\partial_\sigma$. Let us write $z_i=x_i+\sqrt{-1}y_i$. Then we can express the complex structure $J$ on $N$ in terms of the local coordinates as \begin{equation}
	Jdx_i=dy_i,\quad Jdy_i=-dx_i,\quad Jdt=h^{-1}\Theta,
\end{equation} where $\Theta$ is a $1$-form which can be written as \begin{equation}
	\Theta=-d\sigma+\theta
\end{equation} such that $\theta$ does not have $d\sigma$ component, and \begin{equation}
	h^{-1}=Jdt(-\partial_\sigma)=dt(J\partial_\sigma)=|\partial_\sigma|_{\omega_L}^2>0
\end{equation} is independent of $\sigma$. So we write $h=h_t$ as a family of functions varying in $t$.

By the $S^1$-invariance \begin{equation}
	\mathcal{L}_{\partial_\sigma}(Jdt)=0,\quad \mathcal{L}_{\partial_\sigma}d\sigma=0,
\end{equation} we obtain $\mathcal{L}_{\partial_\sigma}\theta=0$. Therefore, $\theta$ is independent of $\sigma$.

Denote $d_X$ the differential along $X$. Since $w$ is holomorphic, $J\partial_{\log|w|}=\xi$, this says \begin{equation}
	dt(\partial_{\log|w|})=Jdt(\xi)=-h_t^{-1}.
\end{equation} By the inverse function theorem, \begin{equation}
	\partial_t\log|w|=-h_t,
\end{equation} hence \begin{equation}
	d\sigma=Jd\log|w|=d_X^c\log|w|+\partial_t\log|w|Jdt=d_X^c\log|w|-\Theta.
\end{equation} So we know $\theta=d_X^c\log|w|$ descends to a family of $1$-forms on $X$ varying in $t$.

We can write the K\"{a}hler form $\omega_L$ on $N$ as \begin{equation}
	\omega=dt\wedge(-d\sigma+\theta)+\omega_t,
\end{equation} where $\omega_t$ is a $(1,1)$-form without $dt$ or $d\sigma$ component. This is due to \labelcref{1} and the fact that $\omega_L$ is of type $(1,1)$. Since $\mathcal{L}_{\partial_\sigma}\omega_L=0$ we also have $\mathcal{L}_{\partial_\sigma}\omega_t=0$, so $\omega_t$ is independent of $\sigma$ and descends to a family of $(1,1)$-forms on $X$ varying in $t$. The condition $d\omega_L=0$ is equivalent to \begin{equation}\label{2}
	\left\{\begin{aligned}
		&d_X\omega_t=0 \\
		&\partial_t\omega_t=d_X\theta.
	\end{aligned}\right.
\end{equation} Therefore, $\omega_t$ is a family of K\"{a}hler metrics on $X$.

We next consider the integrability condition of the complex structure $J$. It is direct to check that \begin{equation}
	J\partial_\sigma=h_t^{-1}\partial_t,
\end{equation} so the holomorphic vector field generating the $\mathbb{C}^*$-action is given by \begin{equation}
	\xi^{1,0}=\frac{1}{2}(\partial_\sigma-\sqrt{-1}J\partial_\sigma)=\frac{1}{2}(\partial_\sigma-\sqrt{-1}h_t^{-1}\partial_t).
\end{equation} The dual holomorphic $(1,0)$-form is \begin{equation}
	\kappa=\sqrt{-1}h_tdt+d\sigma+\kappa',
\end{equation} where $\kappa'$ only involves $dx_i,dy_i$. The integrability condition for $J$ can be expressed as \begin{equation}
	d\kappa\wedge\kappa\wedge dz_1\wedge\cdots\wedge dz_n=0,
\end{equation} which is equivalent to \begin{equation}\label{3}
	\left\{\begin{aligned}
		&d_X\theta\wedge dz_1\wedge\cdots\wedge dz_n=0 \\
		&\partial_t\theta=-d_X^ch_t.
	\end{aligned}\right.
\end{equation} The first equation follows from the second equation in \labelcref{2} which implies $d_X\theta$ is of type $(1,1)$ on $X$. Notice that \labelcref{2,3} together can be re-organized as a system \begin{equation}\label{4}
	\left\{\begin{aligned}
		&\partial_t^2\omega_t+d_Xd_X^ch_t=0 \\
		&d\Theta=\partial_t\omega_t-dt\wedge d_X^ch.
	\end{aligned}\right.
\end{equation}

The $1$-form $-\sqrt{-1}\Theta$ can be viewed as a family of connection $1$-forms on the unit $S^1$-bundle in $L$ over $X$, with curvature \begin{equation}
	-\sqrt{-1}d_X\Theta=-\sqrt{-1}\partial_t\omega_t.
\end{equation} As a consequence $\partial_t\omega_t$ belongs to the cohomology class $2\pi c_1(L)\in 2\pi H^2(X,\mathbb{Z})$.

Now we consider the Calabi--Yau condition. Let $\Omega_X$ be a holomorphic volume form on $X$, then there is a holomorphic volume form \begin{equation}\label{OmegaL}
	\Omega_L=\sqrt{-1}\kappa\wedge\Omega_X
\end{equation} on $L^\times$. The Calabi--Yau condition becomes \begin{equation}\label{5}
	\frac{\omega_L^{n+1}}{c_{n+1}\Omega_L\wedge\overline{\Omega}_L}=\frac{1}{2h_t}\frac{\omega_t^n}{c_n\Omega_X\wedge\overline{\Omega}_X}=1,\quad c_n=n!\sqrt{-1}^{n^2}.
\end{equation}

In summary, \labelcref{4,5} together show the K\"{a}hler family $\omega_t$ is controlled by the generalized Gibbons--Hawking ansatz \begin{equation}\label{gh}
	\left\{\begin{aligned}
		&\partial_t^2\omega_t+\sqrt{-1}\partial_X\bar{\partial}_X(\frac{\omega_t^n}{c_n\Omega_X\wedge\overline{\Omega}_X})=0 \\
		&\partial_t[\omega_t]=2\pi c_1(L).
	\end{aligned}\right.
\end{equation} Once we have determined $\omega_t$, the first equation of \labelcref{4} determines $h_t$ up to a constant varying in $t$. This additive function of $t$ is fixed by imposing \labelcref{5}: \begin{equation}
	\int_X 2h_tc_n\Omega_X\wedge\overline{\Omega}_X=\int_X\omega_t^n
\end{equation}

In particular, we can normalize \begin{equation}
	\int_X c_n\Omega_X\wedge\overline{\Omega}_X=(2\pi c_1(L))^n,
\end{equation} so that \begin{equation}
	\frac{\int_X 2h_tc_n\Omega_X\wedge\overline{\Omega}_X}{\int_X c_n\Omega_X\wedge\overline{\Omega}_X}=t^n+\text{remainder terms}.
\end{equation}

In the case when $L$ is ample, we can find a Calabi--Yau metric $\omega_X\in 2\pi c_1(L)$ on $X$, then $\omega_t=t\omega_X$ solves the equation \labelcref{gh}. And the associated $2h_t=t^n$. This is the classical Calabi model space used in the end of Tian--Yau metric.

If $L$ is not ample, it's hard to solve \labelcref{gh} directly. Alternatively, we construct approximate solutions specified by the data $(\omega_t,h_t)$, with \begin{equation}\label{agh}
	\left\{\begin{aligned}
		&\partial_t^2\omega_t+2\sqrt{-1}\partial_X\bar{\partial}_Xh_t=0 \\
		&\partial_t[\omega_t]=2\pi c_1(L) \\
		&\int_X 2h_tc_n\Omega_X\wedge\overline{\Omega}_X=\int_X\omega_t^n \\
		&\omega_t,h_t>0\text{ for }t\gg1.
	\end{aligned}\right.
\end{equation} The error from being Calabi--Yau is represented by the Ricci potential of $\omega_L$, namely, \begin{equation}
	\log(\frac{1}{2h_t}\frac{\omega_t^n}{c_n\Omega_X\wedge\overline{\Omega}_X}).
\end{equation}

To find the Ricci curvature, we need to compute the global $dd^c$ in terms of $d_Xd_X^c$ and $t$-derivatives. Let $f$ be an $S^1$-invariant function, then \begin{equation}
	\bar{\partial}f=\bar{\partial}_Xf+\frac{1}{2}(dt-\sqrt{-1}Jdt)\partial_tf=\bar{\partial}_Xf+\frac{1}{2}(dt-\sqrt{-1}h_t^{-1}\Theta)\partial_tf.
\end{equation} Denote \begin{equation}
	\eta=\frac{1}{2}(dt+\sqrt{-1}h_t^{-1}\Theta),
\end{equation} then by \labelcref{4} \begin{equation}
	d(h_t\bar{\eta})=\frac{1}{2}(d_Xh_t\wedge dt-\sqrt{-1}d\Theta)=\bar{\partial}_Xh_t\wedge dt-\frac{1}{2}\sqrt{-1}\partial_t\omega_t.
\end{equation} Solve $d\bar{\eta}$ from above and $dd^cf=2\sqrt{-1}d\bar{\partial}f$, we get \begin{equation}\label{ddc}
	dd^cf=d_Xd_X^cf+2\sqrt{-1}h_t(\partial_Xv\wedge\bar{\eta}+\eta\wedge\bar{\partial}_Xv+\partial_tv\eta\wedge\bar{\eta})+v\partial_t\omega_t,
\end{equation} where \begin{equation}
	v=h_t^{-1}\partial_tf.
\end{equation}

In terms of $\eta$, $\omega_L=2\sqrt{-1}h_t\eta\wedge\bar{\eta}+\omega_t$, so we get an estimate \begin{multline}\label{eddc}
	|dd^cf|_{\omega_L}\leq C(|d_Xd_X^cf|_{\omega_t}+h_t^{-1}|\partial_tf|\cdot|\partial_t\omega_t|_{\omega_t} \\
	+h_t^{-1/2}|d_X\partial_tf|_{\omega_t}+h_t^{-3/2}|\partial_tf||d_Xh_t|_{\omega_t}+h_t^{-1}|\partial_t^2f|+h_t^{-2}|\partial_tf||\partial_th_t|).
\end{multline}

\subsection{Recovering the global potential}\label{complement}

We reconstruct the metric on $N$ from the data of (approximate) generalized Gibbons--Hawking ansatz $(\omega_t,h)$.

Choose a Hermitian metric $\|w\|^2=|w|^2e^{-\rho}$ on $L$, $s=-\log\|w\|$, where $w=|w|e^{\sqrt{-1}\sigma}$ is the fiber coordinate. Define \begin{equation}
	\Theta_0:=d^cs=-d\sigma+\frac{1}{2}d_X^c\rho,
\end{equation} and \begin{equation}
	\alpha_X:=d\Theta_0=\sqrt{-1}\partial_X\bar{\partial}_X\rho\in 2\pi c_1(L).
\end{equation}

Since $\partial_t[\omega_t]=2\pi c_1(L)$, we can write $[\omega_t]=2\pi tc_1(L)+[\beta_X]$ for some $\beta_X$ fixed. Hence \begin{equation}
	\omega_t=t\alpha_X+\beta_X+\sqrt{-1}\partial_X\bar{\partial}_X\phi_t
\end{equation} for a family of functions $\phi_t$ on $X$, and we normalize it by \begin{equation}
	\int_X\phi_t\cdot c_n\Omega_X\wedge\overline{\Omega}_X=0.
\end{equation}

From \labelcref{agh} we know \begin{equation}
	2h_t=\frac{(t[\alpha_X]+[\beta_X])^n}{[\alpha_X]^n}-\partial_t^2\phi_t.
\end{equation} We let \begin{equation}
	\Theta=\Theta_0+\frac{1}{2}d_X^c\partial_t\phi_t
\end{equation} so that \begin{equation}
	d\Theta=\partial_t\omega_t-dt\wedge d_X^ch_t.
\end{equation}

Now we identify $s$ with $t$. First notice that they are both $S^1$-invariant. Since \begin{equation}
	2\partial s=ds+\sqrt{-1}\Theta_0=d_Xs+\partial_ts dt+\sqrt{-1}\Theta-\frac{1}{2}\sqrt{-1}d_X^c\partial_t\phi_t
\end{equation} is of type $(1,0)$, \begin{equation}
	\partial_ts dt=-J\Theta=h_tdt.
\end{equation} And therefore, \begin{equation}\label{s}
	s=\frac{1}{2[\alpha_X]^n}\sum_{i=0}^n\binom{n}{i}([\alpha_X]^i[\beta_X]^{n-i})\frac{t^{i+1}}{i+1}-\frac{1}{2}\partial_t\phi_t+c_X,
\end{equation} where $c_X$ is independent of $t$. Substituting this expression into $\partial s$ we find \begin{equation}
	2\partial s=d_Xc_X-\partial_X\partial_t\phi_t+h_tdt+\sqrt{-1}\Theta
\end{equation} is of type $(1,0)$. So $\bar{\partial}_Xc_X=0$ and $c_X=c$ is a constant.

Since $\partial_ts=h_t>0$, $s$ is increasing in $t$. If $s\to+\infty$ as $t\to+\infty$, then the reconstructed metric from $t\gg1$ contains a punctured tubular neighborhood of the zero section. We will verify this in \labelcref{st}.

The reduced data determines the metric as a differential form, but the later
gluing argument requires a potential relative to a fixed representative of
its cohomology class. The following proposition provides this potential
explicitly. In particular, it separates the fixed cohomological contribution
from the part that can be transferred and modified by a
$\sqrt{-1}\partial\bar\partial$-construction.

\begin{proposition}\label{potential}
	The K\"{a}hler metric defined by the (approximate) generalized Gibbons--Hawking ansatz equals \begin{equation}
		\omega_L=dt\wedge\Theta+\omega_t=\pi_L^*\beta_X+\sqrt{-1}\partial\bar{\partial}U,
	\end{equation} where \begin{equation}
		U=2ts-\frac{1}{[\alpha_X]^n}\sum_{i=0}^n\binom{n}{i}([\alpha_X]^i[\beta_X]^{n-i})\frac{t^{i+2}}{(i+1)(i+2)}+\phi_t-2ct.
	\end{equation}
\end{proposition}

\begin{proof}
	Regard $U$ as a function of $t$ and the base coordinates on $X$. By the definition, \begin{equation}
		\partial_t U=2t\partial_t s=2th_t.
	\end{equation} Thus, $h_t^{-1}\partial_tU=2t$. According to \labelcref{ddc}, \begin{equation}
		\begin{aligned}
			dd^cU&=d_Xd_X^cU+4\sqrt{-1}h_t\eta\wedge\bar{\eta}+2t\partial_t\omega_t \\
			&=2\sqrt{-1}\partial_X\bar{\partial}_X(-t\partial_t\phi_t+\phi_t)+2dt\wedge\Theta+2t\alpha_X+2t\sqrt{-1}\partial_X\bar{\partial}_X\partial_t\phi_t \\
			&=2(dt\wedge\Theta+t\alpha_X+\beta_X+\sqrt{-1}\partial_X\bar{\partial}_X\phi_t)-2\beta_X.
		\end{aligned}
	\end{equation}
\end{proof}

\section{The gluing construction}\label{s3}

We now specialize to $X=D$ and construct the initial approximate family for
the reduced system. The orbifold part of the metric expands with $t$,
whereas the areas of the exceptional curves must remain equal to the
prescribed periods $p_i(\mathfrak{c})$. To accommodate these different behaviors, we
first choose cohomological representatives compatible with the ADE
contraction, and then glue the orbifold and ALE models. We introduce
weighted spaces adapted to their different geometric scales and estimate
the error of the resulting family. These estimates provide the starting
point for the correction procedure in \cref{s5}. Throughout this section,
$\partial$ and $\bar\partial$ denote differentiation on $X$.

\subsection{Adjusting the cohomology class}\label{b}

\begin{lemma}
	Let $\mathfrak{c}\in H^2(M,\mathbb{R})$ be a K\"{a}hler class on $M$, such that all exceptional periods $p_i(\mathfrak{c})>0$. Then there exists a K\"{a}hler class $\bar{\mathfrak{c}}\in H^2(\overline{M},\mathbb{R})$, such that $c$ is the restriction of $\bar{\mathfrak{c}}$. That is, let $j\colon M\to\overline{M}$ be the inclusion, then $\mathfrak{c}=j^*\bar{\mathfrak{c}}$.
\end{lemma}

\begin{proof}
	By the Thom--Gysin sequence \begin{equation}
		\begin{tikzcd}
			H^0(D,\mathbb{R}) \ar[r,"{1\mapsto [D]}"] & H^2(\overline{M},\mathbb{R}) \ar[r,"j^*"] & H^2(M,\mathbb{R}) \ar[r] & H^1(D,\mathbb{R})=0,
		\end{tikzcd}
	\end{equation} $\mathfrak{c}$ has a lifting $\bar{\mathfrak{c}}\in H^2(\overline{M},\mathbb{R})$, and all these liftings are represented by $\bar{\mathfrak{c}}+\mathbb{R}[D]$.

	Apply \cite[Proposition 0.6(i)]{demaillyNumericalCharacterizationKahler2004} to the positive current $[D]$, whose positive Lelong locus is $D$. Since $c_1(-K_{\overline{M}})$ restricts to a nef class $c_1(L)$ on $D$, $c_1(-K_{\overline{M}})$ is nef. Since \begin{equation}
		c_1(-K_{\overline{M}})^3=\int_Dc_1(L)^2>0,
	\end{equation} again by \cite[Theorem 0.5]{demaillyNumericalCharacterizationKahler2004} there is a K\"{a}hler current in this class. Then \cite[Remark 3.6]{demaillyNumericalCharacterizationKahler2004} shows $\overline{M}$ is Moishezon. Since $\overline{M}$ is also compact K\"{a}hler, $\overline{M}$ is projective, and hence weak Fano.

	Thus, by Kawamata--Viehweg vanishing, $H^q(\overline{M},\mathcal{O}_{\overline{M}})=0$ for $q\geq1$. Therefore, $\bar{\mathfrak{c}}$ is equivalent to an $\mathbb{R}$-divisor class in $\overline{M}$. By \cite[Theorem 1.4]{gongyoWeakFanoVarieties2012}, \begin{equation}
		\overline{\operatorname{NE}(\overline{M})}_{\mathbb{R}}=\sum_k\mathbb{R}_{\geq0}[C_k]
	\end{equation} is generated by finitely many irreducible curve classes $[C_k]$.

	For any generator $[C_k]$ with $[D]\cdot[C_k]=0$. If $C_k\subset D$, then $(c_1(L)\cdot[C_k])_D=0$. So $C_k=E_i$ is some exceptional component. And hence \begin{equation}
		\bar{\mathfrak{c}}\cdot[C_k]=\bar{\mathfrak{c}}\cdot[E_i]=p_i(\mathfrak{c})>0.
	\end{equation} If $C_k\nsubset D$, since $[D]$ is effective in $\overline{M}$, $C_k\cap D=\emptyset$. So $C_k$ is a compact curve in $M$, which means \begin{equation}
		\bar{\mathfrak{c}}\cdot[C_k]=\int_{C_k}\mathfrak{c}>0.
	\end{equation}

	For other generators, $[D]\cdot[C_k]>0$. So we can choose an $a\gg1$ such that \begin{equation}
		(\bar{\mathfrak{c}}+a[D])\cdot[C_k]>0
	\end{equation} for all generators $[C_k]$. By Kleiman's criterion, $\bar{\mathfrak{c}}+a[D]$ is a K\"{a}hler class that lifts $\mathfrak{c}$.
\end{proof}

From this lemma, we can take a K\"{a}hler form $\kappa$ on $\overline{M}$, such that $[\kappa|_M]=\mathfrak{c}$. Meanwhile, we choose a Hermitian metric on $\mathcal{O}_{\overline{M}}(D)$ extending the chosen one on $L=\mathcal{O}_{\overline{M}}(D)|_D$ in \cref{complement}. Let $\alpha$ be the curvature, so that $\alpha|_D=\alpha_X\in2\pi c_1(L)$. And we write \begin{equation}
	A=\frac{c_1(L)\cdot[\kappa|_D]}{2\pi c_1(L)^2},\quad \beta=\kappa-A\alpha.
\end{equation} Note that $A>0$ from the nefness and bigness assumption. Since the defining section trivializes $\mathcal{O}_{\overline{M}}(D)$ on $M$, $[\alpha|_M]=0$, and $[\beta|_M]=\mathfrak{c}$. The choice of coefficient $A$ makes $c_1(L)\cdot[\beta|_D]=0$. Let $C_{ij}=-E_i\cdot E_j$ be the Cartan matrix, $(C^{ji})=(C_{ij})^{-1}$. Then \begin{equation}\label{class}
	[\beta|_D]=[\gamma]-\sum_ia_i[E_i]:=[\gamma]-\sum_{i,j}p_j(c) C^{ji}[E_i],
\end{equation} where $[\gamma]\cdot c_1(L)=[\gamma]\cdot[E_i]=0$.

\subsection{The ALE and the orbifold models}

Consider $\mathbb{C}^2/\Gamma$ for $\Gamma\subset\operatorname{SU}(2)$ a finite subgroup, and the minimal crepant resolution $Q=\widetilde{\mathbb{C}^2/\Gamma}$. Let $(z_1,z_2)$ be the coordinates on $\mathbb{C}^2$, $r^2=|z_1|^2+|z_2|^2$. $\Omega=dz_1\wedge dz_2$ and lifts to a holomorphic volume form on $Q$. By \cite[Theorem 3.3]{joyceAsymptoticallyLocallyEuclidean2001}, in every prescribed K\"{a}hler class, there exists a unique Ricci-flat ALE K\"{a}hler metric $\omega_Q$ on $Q$, such that $\omega_Q^2=2\Omega\wedge\overline{\Omega}$, and for $r\geq C$ on the end, \begin{equation}
	\omega_Q=\sqrt{-1}\partial\bar{\partial}\Psi_Q,\quad|\nabla^k(\Psi_Q-|z|^2)|\leq C|z|^{-2-k}.
\end{equation}

We can require the potential $\Psi_Q$ to be circle invariant. Let $R=J(r\partial_r)$ the Reeb field with respect to the flat cone metric on $\mathbb{C}^2/\Gamma$. In coordinates, \begin{equation}
	R=2\operatorname{Re}(\sqrt{-1}(z_1\frac{\partial}{\partial z_1}+z_2\frac{\partial}{\partial z_2})).
\end{equation} So $R$ is holomorphic and the flow of $R$ is the circle action \begin{equation}
	R_\theta(z)=e^{\sqrt{-1}\theta}z.
\end{equation}

Since \begin{equation}
	\mathcal{L}_R\omega_Q=d\iota_R\omega_Q+\iota_Rd\omega_Q=d\iota_R\omega_Q,
\end{equation} $R$ preserves the K\"{a}hler class. Since $\mathcal{L}_Rr=0$, $R$ also preserves the asymptote of $\omega_Q$. The uniqueness part of the theorem shows $R_\theta^*\omega_Q=\omega_Q$.

Therefore, on the end, by averaging along the circle direction, we let \begin{equation}
	\widetilde{\Psi}_Q=\frac{1}{2\pi}\int_0^{2\pi}R_\theta^*\Psi_Qd\theta.
\end{equation} It follows that \begin{equation}
	\sqrt{-1}\partial\bar{\partial}\widetilde{\Psi}_Q=\frac{1}{2\pi}\int_0^{2\pi}R_\theta^*\sqrt{-1}\partial\bar{\partial}\Psi_Qd\theta=\frac{1}{2\pi}\int_0^{2\pi}R_\theta^*\omega_Qd\theta=\omega_Q.
\end{equation}

From now on, we always assume the end potential $\Psi_Q$ is circle invariant.

\begin{proposition}\label{psi}
	The circle invariant $\Psi_Q$ extends smoothly to $\{r>0\}$, that is, the complement of the exceptional locus in $Q$, keeping $\omega_Q=\sqrt{-1}\partial\bar{\partial}\Psi_Q$.
\end{proposition}

\begin{proof}
	First we extend the radial derivative $r\partial_r\Psi_Q$. On the end, since $\Psi_Q$ is circle invariant, $R\Psi_Q=0$. So \begin{equation}
		\begin{aligned}
			d(r\partial_r\Psi_Q)&=-d(\iota_{JR}d\Psi_Q) \\
			&=d\iota_RJd\Psi_Q \\
			&=\mathcal{L}_Rd^c\Psi_Q-\iota_Rdd^c\Psi_Q \\
			&=d^c(R\Psi_Q)-2\iota_R\sqrt{-1}\partial\bar{\partial}\Psi_Q \\
			&=-2\iota_R\omega_Q.
		\end{aligned}
	\end{equation}

	From previous analysis, $R_\theta$ preserves $\omega_Q$, hence $\mathcal{L}_R\omega_Q=d\iota_R\omega_Q=0$. We can show there is a global moment function. In fact, by the ALE Hodge theory, \cite[Theorem 5.2]{joyceAsymptoticallyLocallyEuclidean2001}, we can write \begin{equation}
		\iota_R\omega_Q=\eta+df
	\end{equation} for a unique harmonic representative $\eta$, satisfying \begin{equation}
		d\eta=0,\quad d^*\eta=0,\quad|\eta|\leq Cr^{-3}.
	\end{equation} Since $\omega_Q$ is Ricci-flat, the Bochner formula gives \begin{equation}
		\Delta|\eta|^2=|\nabla\eta|^2\geq0.
	\end{equation} The maximum principle, together with the decaying property, shows that $\eta=0$. So $\iota_R\omega_Q=df$, where $f$ is globally defined on $Q$. After fixing the additive constant we can extend \begin{equation}
		r\partial_r\Psi_Q=-2f.
	\end{equation}

	The flow of $r\partial_r$ is the radial dilation \begin{equation}
		\delta_s(z)=e^sz.
	\end{equation} Thus, we know \begin{equation}
		\frac{d}{ds}\delta_s^*\Psi_Q=\delta_s^*(r\partial_r\Psi_Q)=-2\delta_s^*f.
	\end{equation} For any point not in the exceptional locus, we can choose $A>0$ large so that $e^Ar\geq C$, then \begin{equation}
		\Psi_Q=\delta_A^*\Psi_Q+2\int_0^A\delta_s^*fds
	\end{equation} extends the definition of $\Psi_Q$. To see $\omega_Q=\sqrt{-1}\partial\bar{\partial}\Psi_Q$, first we compute \begin{equation}
		dd^cf=dJ\iota_R\omega_Q=d\iota_{-r\partial_r}\omega_Q=-\mathcal{L}_{r\partial_r}\omega_Q.
	\end{equation} Therefore, \begin{equation}
		\begin{aligned}
			\sqrt{-1}\partial\bar{\partial}\Psi_Q&=\delta_A^*\sqrt{-1}\partial\bar{\partial}\Psi_Q+\int_0^A\delta_s^*(dd^cf)ds \\
			&=\delta_A^*\omega_Q-\int_0^A\delta_s^*(\mathcal{L}_{r\partial_r}\omega_Q)ds \\
			&=\delta_A^*\omega_Q-\int_0^A\frac{d}{ds}(\delta_s^*\omega_Q)ds \\
			&=\omega_Q.
		\end{aligned}
	\end{equation}
\end{proof}

Since $\delta_s^*\Omega=e^{2s}\Omega$, we rescale the metric by \begin{equation}
	\omega_{Q,a}=a\delta_{-\frac{1}{2}\log a}^*\omega_Q,\quad\Psi_{Q,a}=a\delta_{-\frac{1}{2}\log a}^*\Psi_Q.
\end{equation} Then the Calabi--Yau relation \begin{equation}
	\omega_{Q,a}^2=2\Omega\wedge\overline{\Omega}
\end{equation} holds. Because $\delta_s$ fixes the exceptional locus, \begin{equation}
	\int_{E_i}\omega_{Q,a}=\int_{E_i}a\delta_{-\frac{1}{2}\log a}^*\omega_Q=\int_{E_i}a\omega_Q,
\end{equation} the K\"{a}hler class is rescaled by \begin{equation}
	[\omega_{Q,a}]=a[\omega_Q].
\end{equation}

Normalize the holomorphic volume form $\Omega_X$ on $X$ as \begin{equation}
	\int_X 2\Omega_X\wedge\overline{\Omega}_X=4\pi^2c_1(L)^2,
\end{equation} and let $\Omega_Y$ denote its descent on $Y$. Since $c_1(L)$ descends to a positive orbifold $(1,1)$-class on $Y$. From the orbifold Calabi--Yau theorem \cite[Theorem 1.2]{faulkYausTheoremEffective2019}, there exists an orbifold Ricci-flat K\"{a}hler metric $\omega_Y\in2\pi c_1(L)$, \begin{equation}
	\omega_Y^2=2\Omega_Y\wedge\overline{\Omega}_Y.
\end{equation}

Also, the previous class $[\gamma]$ descends to an orbifold $(1,1)$-class on $Y$. Again from the orbifold Calabi--Yau theorem, for $a>0$ small enough such that the class $2\pi c_1(L)+a[\gamma]$ is positive, there exists a family of orbifold Ricci-flat K\"{a}hler metrics $\omega_{Y,a}\in2\pi c_1(L)+a[\gamma]$, \begin{equation}
	\omega_{Y,a}^2=(1+\frac{[\gamma]^2}{4\pi^2 c_1(L)^2}a^2)\cdot2\Omega_Y\wedge\overline{\Omega}_Y.
\end{equation} We can write \begin{equation}
	\omega_{Y,a}=\omega_Y+a\gamma_Y+\sqrt{-1}\partial\bar{\partial}u_a,\quad\int_Xu_a\cdot2\Omega_Y\wedge\overline{\Omega}_Y=0,
\end{equation} where $\gamma_Y$ is an orbifold $(1,1)$-form in $[\gamma]$. The implicit function theorem shows that $\omega_{Y,a}$ varies smoothly in $a$. Taking the Taylor expansion, and absorb the linear term by changing $\gamma_Y$, we can require \begin{equation}
	u_a=O(a^2).
\end{equation}

\subsection{The gluing argument}\label{3.2}

Let $(p_i,U_i)$ be the orbifold chart near each singularity on $Y$, $U_i=B/\Gamma_i$, $B\subset\mathbb{C}^2$. Take coordinates $(z_{1,i},z_{2,i})$ on $B$, such that $\Omega_Y=dz_{1,i}\wedge dz_{2,i}$, $r_i^2=|z_{1,i}|^2+|z_{2,i}|^2$. Without loss of generality, if $E_1,\ldots,E_k$ are the exceptional components from resolving $p_i$. Then on the model space $Q_i$, \begin{equation}
	-\sum_{j=1}^ka_j[E_j]
\end{equation} is a K\"{a}hler class, where the coefficients are given by \labelcref{class}. Let $\omega_{Q_i}$ be the associated ALE K\"{a}hler metric in this class, and the rescaled metric $\omega_{Q_i,a}$. Meanwhile, we can adjust the coordinates so that the orbifold metric near $p_i$ has \begin{equation}
	\omega_Y=\sqrt{-1}\partial\bar{\partial}(r_i^2+u_{Y,i}),\quad\gamma_Y=\sqrt{-1}\partial\bar{\partial}u_{\gamma,i},
\end{equation} where $u_{Y,i}=O(r_i^4)$, $u_{\gamma,i}=O(r_i^2)$ near $p_i$.

Repeat this process for each singularity, and we take the orbifold charts $U_i$ disjoint.

Now we take a cut-off function \begin{equation}
	\chi_0=\left\{\begin{aligned}
		&0,&&(-\infty,1] \\
		&1,&&[2,+\infty),
	\end{aligned}\right.
\end{equation} and let $\chi_i=\chi_0(a^{-1/3}r_i)$. We define \begin{equation}
	\hat{\omega}_a:=\left\{\begin{aligned}
		&\sqrt{-1}\partial\bar{\partial}(\chi_i(r_i^2+u_{Y,i}+au_{\gamma,i}+u_a)+(1-\chi_i)\Psi_{Q_i,a}),&&\text{in each }\pi_Y^{-1}(U_i) \\
		&\pi_Y^*\omega_{Y,a},&&\text{elsewhere},
	\end{aligned}\right.
\end{equation} where $\sqrt{-1}\partial\bar{\partial}\Psi_{Q_i,a}$ extends smoothly to the exceptional locus.

If there's no ambiguity, we drop the subscript $i$ and let \begin{equation}
	u_{Y,a}=u_Y+au_\gamma+u_a.
\end{equation} The potential \begin{equation}
	\chi(r^2+u_{Y,a})+(1-\chi)\Psi_{Q,a}=r^2+\chi u_{Y,a}+(1-\chi)(\Psi_{Q,a}-r^2),
\end{equation} where \begin{equation}
	u_{Y,a}=O(r^4)+aO(r^2)+O(a^2)
\end{equation} and \begin{equation}
	\Psi_{Q,a}-r^2=a\delta_{-\frac{1}{2}\log a}^*(\Psi_Q-r^2)=O(\frac{a^2}{r^2}).
\end{equation} Therefore, in the gluing region $\{a^{1/3}\leq r\leq2a^{1/3}\}$, relative to the Euclidean metric $\omega_{\mathrm{Euc}}=\sqrt{-1}\partial\bar{\partial}r^2$, \begin{equation}
	r^k|\nabla^k(\hat{\omega}_a-\omega_{\mathrm{Euc}})|_{\omega_{\mathrm{Euc}}}\leq C a^{2/3}.
\end{equation} This estimate also shows $\hat{\omega}_a$ is K\"{a}hler for $a\ll1$ small enough.

Now we compute the cohomology class of $\hat{\omega}_a$. Since \begin{equation}
	\hat{\omega}_a-\pi_Y^*\omega_{Y,a}=\sqrt{-1}\partial\bar{\partial}\left(\sum_i(1-\chi_i)(\Psi_{Q_i,a}-r_i^2-u_{Y,i}-au_{\gamma,i}-u_a)\right),
\end{equation} with a smooth potential on $X\setminus\bigcup_i\pi_Y^{-1}(p_i)$, we find \begin{equation}
	[\hat{\omega}_a]-[\pi_Y^*\omega_{Y,a}]\in\ker(H^2(X,\mathbb{R})\to H^2(X\setminus\bigcup_jE_j,\mathbb{R})).
\end{equation} By the Thom--Gysin sequence \begin{equation}
	\begin{tikzcd}
		\displaystyle\bigoplus_j\mathbb{R}[E_j] \ar[r] & H^2(X,\mathbb{R}) \ar[r] & H^2(X\setminus\displaystyle\bigcup_jE_j,\mathbb{R}),
	\end{tikzcd}
\end{equation} $[\hat{\omega}_a]-[\pi^*\omega_{Y,a}]$ is a linear combination of $[E_j]$. From the Cartan matrix, we get \begin{equation}
	[\hat{\omega}_a]=[\pi^*\omega_{Y,a}]-\sum_{i,j}(\int_{E_j}\hat{\omega}_a)C^{ji}[E_i]=2\pi c_1(L)+a[\gamma]-a\sum_{i,j}p_j(c)C^{ji}[E_i]=2\pi c_1(L)+a[\beta|_X].
\end{equation}

Let $a=1/t$ and $\omega_t=t\hat{\omega}_{1/t}$, then \begin{equation}
	[\omega_t]=2\pi tc_1(L)+[\beta|_X],
\end{equation} and $\omega_t$ is K\"{a}hler for $t\gg1$ large enough.

We can write $\omega_t$ in terms of a global potential. To see this, fix a constant $C_0>0$ such that the cut-off function \begin{equation}
	\tilde{\chi}_i=\chi_0(C_0r_i)=\left\{\begin{aligned}
		&0,&&\text{near }\pi_Y^{-1}(p_i) \\
		&1,&&\text{outside }\pi_Y^{-1}(U_i).
	\end{aligned}\right.
\end{equation} Thus, for $t\gg1$, the gluing region $r_i\sim t^{-1/3}$ is contained in $\tilde{\chi}_i=0$. Let \begin{equation}\label{13}
	\alpha_0=\left\{\begin{aligned}
		&\sqrt{-1}\partial\bar{\partial}(\tilde{\chi}_i(r_i^2+u_{Y,i})),&&\text{in each }\pi_Y^{-1}(U_i) \\
		&\pi_Y^*\omega_Y,&&\text{elsewhere},
	\end{aligned}\right.
\end{equation} and \begin{equation}\label{18}
	\beta_0=\left\{\begin{aligned}
		&\sqrt{-1}\partial\bar{\partial}(\tilde{\chi}_iu_{\gamma,i}+(1-\tilde{\chi}_i)\Psi_{Q_i}),&&\text{in each }\pi_Y^{-1}(U_i) \\
		&\pi_Y^*\gamma_Y,&&\text{elsewhere}.
	\end{aligned}\right.
\end{equation} Similarly, from the Thom--Gysin sequence we find \begin{equation}
	[\alpha_0]=2\pi c_1(L),\quad[\beta_0]=[\beta|_X].
\end{equation} Therefore, $\omega_t=t\alpha_0+\beta_0+\sqrt{-1}\partial\bar{\partial}v_t$, where \begin{multline}
	v_t=\sum_i(1-\tilde{\chi}_i)(t(r_i^2+u_{Y,i})+u_{\gamma,i}-\Psi_{Q_i}) \\
	+\sum_i(1-\chi_i)(\delta_{\frac{1}{2}\log t}^*\Psi_{Q_i}-t(r_i^2+u_{Y,i})-u_{\gamma,i})+tu_{1/t}-\sum_i(1-\chi_i)tu_{1/t}.
\end{multline} is a smooth global potential. For later use, we give a rough estimate of $v_t$.

\begin{proposition}\label{phi}
	We have \begin{equation}
		t^k\|\partial_t^k v_t\|_{C^0}\leq C_kt
	\end{equation} for $k\geq0$.
\end{proposition}

\begin{proof}
	We estimate separately in four regions.

	Region 1: $\{r<t^{-1/3}\}$.

	By assumption \begin{equation}
		v_t=\delta_{\frac{1}{2}\log t}^*\Psi_Q-\Psi_Q=\int_0^{\frac{1}{2}\log t}\delta_s^*(r\partial_r\Psi_Q)ds.
	\end{equation} Since $r\partial_r\Psi_Q$ is smooth in the core, and of $O(r^2)$ as $r$ tends to infinity, we have \begin{equation}
		|r\partial_r\Psi_Q|\leq C(1+r^2),
	\end{equation} and that \begin{equation}
		|\delta_{\frac{1}{2}\log t}^*\Psi_Q-\Psi_Q|\leq C\int_0^{\frac{1}{2}\log t}(1+e^{2s}r^2)ds\leq C(\log t+tr^2)\leq C(\log t+t^{1/3})
	\end{equation}

	Region 2: $\{t^{-1/3}\leq r\leq2t^{-1/3}\}$.

	\begin{equation}
		\begin{aligned}
			v_t&=\delta_{\frac{1}{2}\log t}^*\Psi_Q-\Psi_Q-\chi(\delta_{\frac{1}{2}\log t}^*\Psi_Q-t(r^2+u_Y)-u_\gamma-tu_{1/t}) \\
			&=\delta_{\frac{1}{2}\log t}^*\Psi_Q-\Psi_Q-\chi(\delta_{\frac{1}{2}\log t}^*(\Psi_Q-r^2)-tu_Y-u_\gamma-tu_{1/t}),
		\end{aligned}
	\end{equation} where \begin{align}
		&\delta_{\frac{1}{2}\log t}^*\Psi_Q-\Psi_Q=O(\log t+tr^2)=O(t^{1/3}), \\
		&\delta_{\frac{1}{2}\log t}^*(\Psi_Q-r^2)=O(t^{-1}r^{-2})=O(t^{-1/3}), \\
		&tu_Y=O(tr^4)=O(t^{-1/3}), \\
		&u_\gamma=O(r^2)=O(t^{-2/3}), \\
		&tu_{1/t}=O(t^{-1}).
	\end{align}

	Region 3: $\{2t^{-1/3}<r<2C_0^{-1}\}$.

	\begin{equation}
		\begin{aligned}
			v_t&=(1-\tilde{\chi})(t(r^2+u_Y)+u_\gamma-\Psi_Q)+tu_{1/t}  \\
			&=(1-\tilde{\chi})(tu_Y+u_\gamma-\delta_{\frac{1}{2}\log t}^*(\Psi_Q-r^2)+(\delta_{\frac{1}{2}\log t}^*\Psi_Q-\Psi_Q))+tu_{1/t},
		\end{aligned}
	\end{equation} where \begin{align}
		&\delta_{\frac{1}{2}\log t}^*\Psi_Q-\Psi_Q=O(\log t+tr^2)=O(t), \\
		&\delta_{\frac{1}{2}\log t}^*(\Psi_Q-r^2)=O(t^{-1}r^{-2})=O(t^{-1/3}), \\
		&tu_Y=O(tr^4)=O(t), \\
		&u_\gamma=O(r^2)=O(1), \\
		&tu_{1/t}=O(t^{-1}).
	\end{align}

	Region 4: Away from the singular charts, that is, every $\tilde{\chi}_i=1$.

	In this region, $v_t=tu_{1/t}=O(t^{-1})$.

	These estimates improve as taking derivatives, so the result.
\end{proof}

\subsection{Weighted spaces adapted to $\omega_t$}\label{3.3}

We use $g_t,\nabla_t,\dist_t$ to denote the Riemannian metric, the Levi-Civita connection, and the distance function with respect to $\omega_t$.

In $\{r<t^{-1/3}\}$, $\omega_t=\delta_{\frac{1}{2}\log t}^*\omega_Q$. As $t\to+\infty$, the gluing scale $r\sim t^{-1/3}$ lies in the asymptotic region of the ALE K\"{a}hler metric: $tr^2\sim t^{1/3}\to+\infty$. We use Joyce's radius function \cite{joyceAsymptoticallyLocallyEuclidean2001} to capture the ALE behavior. That is, a function valued in $[1,\infty)$ and comparable to $\dist_t(\cdot,\bigcup_jE_j)\sim\sqrt{t}r$ if $tr^2\gg 1$. Outside the gluing region, the ordinary norm on $(Y,\omega_Y)$ scaled by metric dilation $\omega_t\sim t\pi^*\omega_Y$ is considered.

We let $\tilde{r}$ be smooth away from exceptional locus, that agrees with $r_i$ in each ALE core, and with the scale $\tilde{r}\sim1$ elsewhere. Then we define the weight function \begin{equation}
	l_t:=1+\sqrt{t}\tilde{r}\sim\left\{\begin{aligned}
		&1,&&r_i\lesssim t^{-1/2} \\
		&\sqrt{t}r_i,&&t^{-1/2}\lesssim r_i\lesssim1 \\
		&\sqrt{t},&&\tilde{r}\sim1.
	\end{aligned}\right..
\end{equation}

The weight $l_t$ records the local geometric scale of the family: it remains
bounded in the ALE cores and is comparable to $\sqrt{t}$ in the outer
orbifold region. The next lemma shows that rescaling by this local length
scale gives uniformly controlled geometry. This is the geometric input
for local elliptic estimates with constants independent of $t$, and
motivates the weighted H\"older norms introduced below.

\begin{lemma}\label{qatlas}
	For a small constant $c>0$, the balls $B=(B_{l_t(x)^{-2}g_t}(x,c),l_t(x)^{-2}g_t)$ has a uniform geometry up to every finite order for every $x\in X$ and $t$.
\end{lemma}

\begin{proof}
	We prove by regions. Fix a large constant $K>0$.

	Region 1: $r(x)\leq Kt^{-1/2}$.

	Use the coordinate transform \begin{equation}
		\tau_t\colon\{|y|<t^{1/6}\}\to\{r<t^{-1/3}\},\quad y\mapsto z=t^{-1/2}y,
	\end{equation} then $\tau_t^*\omega_t=\omega_Q$. Let $x'=\tau_t^{-1}(x)$. Hence the ball \begin{equation}
		B\sim(B(x',c),g_Q).
	\end{equation}

	Region 2: $Kt^{-1/2}\leq r(x)\leq K^{-1}$.

	Consider \begin{equation}
		\tilde{\tau}_t\colon\zeta\mapsto z=z(x)+r(x)\zeta,
	\end{equation} and $x''=\tilde{\tau}_t^{-1}(x)$. We've known \begin{equation}
		r^k|\nabla^k(\frac{1}{t}\omega_t-\sqrt{-1}\partial\bar{\partial}r^2)|_{\sqrt{-1}\partial\bar{\partial}r^2}\leq C(r^2+\frac{1}{t^2r^4}+\frac{1}{t}+\frac{1}{t^2r^2})\leq \frac{C}{K^2}.
	\end{equation} Therefore, by choose the $K$ large enough, \begin{equation}
		l_t(x)^{-2}\tilde{\tau}_t^*\omega_t\sim\sqrt{-1}\partial\bar{\partial}|\zeta|^2=\omega_{\mathrm{Euc}}.
	\end{equation} And \begin{equation}
		B\sim(B(x'',c),g_{\mathrm{Euc}}).
	\end{equation}

	Region 3: $\tilde{r}(x)>K^{-1}$.

	\begin{equation}
		B\sim(B(x,c),g_{Y,1/t}),
	\end{equation} which converges smoothly to $g_Y$ away from the singularities.

	In summary we get the uniformly bounded geometry.
\end{proof}

The uniform geometry in every order also implies the injectivity radius $\operatorname{injrad}(x,l_t(x)^{-2}g_t)\geq i_0$ is uniformly bounded below. Hence we can define the weighted H\"{o}lder norm \begin{equation}
	\|u\|_{C_{\delta,t}^{k,\alpha}}:=\sum_{i=0}^k\sup_Xl_t^{i-\delta}|\nabla_t^iu|_{g_t}+\sup_{0<\dist_t(x,y)\leq\frac{i_0}{2}l_t(x)}l_t^{k+\alpha-\delta}\frac{|\nabla_t^ku(x)-P_{y,x}^t\nabla_t^ku(y)|_{g_t}}{\dist_t(x,y)^\alpha}
\end{equation} for a function $u\colon X\to\mathbb{R}$, where $P_{y,x}^t$ is the parallel transport from $y$ to $x$, along the unique shortest geodesic. If $u$ is varying in $t$, we define \begin{equation}
	\|u\|_{C_{\delta,t}^{k,\alpha;m}}:=\sum_{i=0}^mt^i\|\partial_t^iu\|_{C_{\delta,t}^{k,\alpha}}.
\end{equation}

Denote \begin{equation}
	P_t=\chi(tr^2+tu_Y+u_\gamma+tu_{1/t})+(1-\chi)\delta_{\frac{1}{2}\log t}^*\Psi_Q
\end{equation} near each singularity. Using the notation in \cref{psi}, we know \begin{equation}
	\partial_t^2\delta_{\frac{1}{2}\log t}^*\Psi_Q=\frac{1}{4t^2}\delta_{\frac{1}{2}\log t}^*((r\partial_r)^2\Psi_Q-2r\partial_r\Psi_Q).
\end{equation} So $\partial_t^2P_t$ extends smoothly to $X$, and \begin{equation}
	\partial_t^2\omega_t=\sqrt{-1}\partial\bar{\partial}\partial_t^2P_t.
\end{equation} From \labelcref{agh}, the associated $h_t$ to $\omega_t$ is \begin{equation}
	2h_t=t^2+\frac{[\beta|_X]^2}{4\pi^2c_1(L)^2}-\partial_t^2P_t+\operatorname{Avg}(\partial_t^2P_t),
\end{equation} where \begin{equation}
	\operatorname{Avg}(u):=\frac{\int_X u\cdot2\Omega\wedge\overline{\Omega}}{\int_X 2\Omega\wedge\overline{\Omega}}.
\end{equation} Next, we'll estimate the error of $(\omega_t,h_t)$ of being Calabi--Yau, in terms of the weighted norms.

\begin{proposition}\label{F}
	Let $F_t=\log(\frac{\omega_t^2}{t^2\cdot2\Omega\wedge\overline{\Omega}})$, then there is \begin{equation}
		\|F_t\|_{C_{0,t}^{k,\alpha;m}}\leq C_{k,\alpha,m}t^{-2/3}
	\end{equation} for every $k,m\geq0$, $0<\alpha<1$.
\end{proposition}

\begin{proof}
	Region 1: $\{\tilde{r}>2t^{-1/3}\}$.

	\begin{equation}
		F_t=\log(1+\frac{[\gamma]^2}{4\pi^2c_1(L)^2}\frac{1}{t^2}),
	\end{equation} the estimate holds clearly here.

	Region 2: $\{t^{-1/3}\leq r\leq2t^{-1/3}\}$.

	Use the coordinate change \begin{equation}
		\sigma_t\colon\{1\leq|w|\leq2\}\to\{t^{-1/3}\leq r\leq2t^{-1/3}\},\quad w\mapsto z=t^{-1/3}w,
	\end{equation} we know \begin{equation}
		t^{-1/3}\sigma_t^*\omega_t=\sqrt{-1}\partial\bar{\partial}(t^{-1/3}P_t(t^{-1/3}w)).
	\end{equation}

	Let $\epsilon=t^{-1/3}$ and \begin{multline}
		Q_\epsilon(w):=t^{-1/3}P_t(t^{-1/3}w) \\
		=|w|^2+\chi_0(|w|)(\epsilon^{-2}u_Y(\epsilon w)+\epsilon u_\gamma(\epsilon w)+\epsilon^{-2} u_{\epsilon^3}(\epsilon w))+(1-\chi_0(|w|))(\epsilon\Psi_Q(\epsilon^{-1/2}w)-|w|^2).
	\end{multline}

	By Taylor expansion, we can write \begin{equation}
		u_Y(\epsilon w)=\frac{\epsilon^4}{6}\int_0^1(1-y)^3u_Y^{(4)}(\epsilon yw)w^4dy.
	\end{equation} Similarly, $u_\gamma(\epsilon w)$ and $u_{\epsilon^3}(\epsilon w)$ have orders $\epsilon^2$ and $\epsilon^6$, respectively. Meanwhile, we have \begin{equation}
		(\epsilon\partial_\epsilon)^m(\epsilon\Psi_Q(\epsilon^{-1/2}w)-|w|^2)=\epsilon((1-\frac{1}{2}r\partial_r)^m(\Psi_Q-r^2))(\epsilon^{-1/2}w),
	\end{equation} hence \begin{equation}
		\|(\epsilon\partial_\epsilon)^m(\epsilon\Psi_Q(\epsilon^{-1/2}w)-|w|^2)\|_{C^k(\{1\leq|w|\leq2\})}\leq C_{m,k}\epsilon^2
	\end{equation} In summary, \begin{equation}
		Q_\epsilon(w)=|w|^2+\epsilon^2 R(\epsilon,w),
	\end{equation} where $R$ has uniformly controlled derivatives in $(\log\epsilon,w)$. Thus, \begin{equation}\label{17}
		t^{-1/3}\sigma_t^*\omega_t=\omega_{\mathrm{Euc}}+\epsilon^2\sqrt{-1}\partial\bar{\partial}R(\epsilon,w),
	\end{equation} and \begin{equation}
		F_t=\log\det(1+\epsilon^2\partial_{i\bar{j}}R)=\epsilon^2\int_0^1\operatorname{tr}((1+y\epsilon^2\partial_{i\bar{j}}R)^{-1}\partial_{i\bar{j}}R)dy=:\epsilon^2G.
	\end{equation} $G$ is also uniformly controlled in derivatives in $(\log\epsilon,w)$. Applying induction to \begin{equation}
		t\partial_t G=-\frac{1}{3}\epsilon\partial_\epsilon G+\frac{1}{3}w\cdot\partial_wG,
	\end{equation} we show that \begin{equation}
		t^m\|\partial_t^mG\|_{C^k(\{1\leq|w|\leq2\})}\leq C_{m,k}.
	\end{equation} Back to the original coordinates, recall $l_t\sim t^{1/6}$ as $r\sim t^{-1/3}$, and bound the H\"{o}lder term by higher-ordered differential, hence the conclusion.
\end{proof}

\Cref{F} controls the volume form error of $\omega_t$ relative to its
leading $t^2$-normalization. The reconstructed metric, however, involves the
factor $2h_t$, so we must also compare $2h_t$ with $t^2$. \Cref{h} supplies this estimate. \Cref{coe} then controls the
parameter derivatives of the metric, which will be needed when
differentiating the Poisson equation in \cref{s5}.

\begin{proposition}\label{h}
	There is \begin{equation}
		\|\frac{2h_t}{t^2}-1\|_{C_{0,t}^{k,\alpha;m}}\leq C_{k,\alpha,m}t^{-2}
	\end{equation} for every $k,m\geq0$, $0<\alpha<1$.
\end{proposition}

\begin{proof}
	In gluing region $\{t^{-1/3}\leq r\leq2t^{-1/3}\}$, we've known \begin{equation}
		P_t=tr^2+\epsilon R(\epsilon,w),
	\end{equation} hence \begin{equation}
		\partial_t^2P_t=\partial_t^2(\epsilon R(\epsilon,w))=O(t^{-7/3}).
	\end{equation} Here the big $O$ notation improves in the weighted norm.

	In $\{r<t^{-1/3}\}$, \begin{equation}
		\partial_t^2P_t=\frac{1}{4t^2}\delta_{\frac{1}{2}\log t}^*((r\partial_r)^2\Psi_Q-2r\partial_r\Psi_Q)=O(t^{-2}),
	\end{equation} and \begin{equation}
		\delta_{\frac{1}{2}\log t}^*((r\partial_r)^2\Psi_Q-2r\partial_r\Psi_Q)\leq C\min(1,\frac{1}{tr^2})\leq Cl_t^{-2}.
	\end{equation} So we calculate \begin{equation}
		\left|\int_{r<t^{-1/3}}\partial_t^2P_t\cdot2\Omega\wedge\overline{\Omega}\right|\leq C\int_0^{t^{-1/3}}\frac{r^3dr}{t^2(1+\sqrt{t}r)^2}=Ct^{-4}\int_0^{t^{1/6}}\frac{y^3dy}{(1+y)^2}\leq Ct^{-11/3},
	\end{equation} and \begin{equation}
		\left|\int_{t^{-1/3}\leq r\leq 2t^{-1/3}}\partial_t^2P_t\cdot2\Omega\wedge\overline{\Omega}\right|\leq Ct^{-7/3}\int_{t^{-1/3}}^{2t^{-1/3}}r^3dr\leq Ct^{-11/3}.
	\end{equation}

	In $\{\tilde{r}>2t^{-1/3}\}$, only the term $tu_{1/t}$ survives after taking twice differentiation in $t$, with an order $O(t^{-3})$. And these estimates improve as taking derivative in $t$, so we get \begin{equation}
		t^m|\partial_t^m\operatorname{Avg}(\partial_t^2P_t)|\leq C_mt^{-3}.
	\end{equation}

	Therefore, \begin{equation}
		\begin{aligned}
			\|\frac{2h_t}{t^2}-1\|_{C_{0,t}^{k,\alpha;m}}&\leq\|\frac{[\beta|_X]^2}{4\pi^2c_1(L)^2}\frac{1}{t^2}\|_{C_{0,t}^{k,\alpha;m}}+\|\frac{\partial_t^2P_t}{t^2}\|_{C_{0,t}^{k,\alpha;m}}+\|\frac{\operatorname{Avg}(\partial_t^2P_t)}{t^2}\|_{C_{0,t}^{k,\alpha;m}} \\
			&\leq C_{k,\alpha,m}(t^{-2}+t^{-4}+t^{-5})\leq C_{k,\alpha,m}t^{-2},
		\end{aligned}
	\end{equation}
\end{proof}

\begin{proposition}\label{coe}
	There is \begin{equation}
		l_t^it^j|\nabla_{\omega_t}^i\partial_t^j\omega_t|_{\omega_t}\leq C_{i,j}
	\end{equation} for every $i,j\geq0$.
\end{proposition}

\begin{proof}
	 The proof again proceeds by considering the three regions, we briefly sketch it.

	 In $\{r<t^{-1/3}\}$, $\omega_t=\delta_{\frac{1}{2}\log t}^*\omega_Q$, \begin{equation}
	 	l_t^i\nabla_{\omega_t}^i(t\partial_t)^j\omega_t=\delta_{\frac{1}{2}\log t}^*\left((1+r)^i\nabla_{\omega_Q}^i(\frac{1}{2}\mathcal{L}_{r\partial_r})^j\omega_Q\right).
	 \end{equation} Then induce on $j$.

	 In $\{t^{-1/3}\leq r\leq2t^{-1/3}\}$, use \labelcref{17}, the uniform estimate of $R(\epsilon,w)$, and the derivative formula \begin{equation}
	 	t\partial_t R=-\frac{1}{3}\epsilon\partial_\epsilon R+\frac{1}{3}w\cdot\partial_wR.
	 \end{equation}

	 In $\{\tilde{r}>2t^{-1/3}\}$, $\omega_t=t\omega_{Y,1/t}$, $l_t\leq C\sqrt{t}$. So \begin{equation}
	 	t^{i/2}t^j|\nabla_{\omega_t}^i\partial_t^j\omega_t|_{\omega_t}=t^{j-1}|\nabla_{\omega_{Y,1/t}}^i(t\partial_t^j\omega_{Y,1/t}+j\partial_t^{j-1}\omega_{Y,1/t})|_{\omega_{Y,1/t}}.
	 \end{equation} Then use $t\partial_t\omega_{Y,1/t}=-a\partial_a\omega_{Y,a}$ and the smoothness of $\omega_{Y,a}$ at $a=0$.
\end{proof}

\section{Inverting the Laplacian}\label{s4}

This section supplies the linear analytic tool needed to improve the
approximate solution. Solvability of the Poisson equation on each fixed
compact manifold $(X,\omega_t)$ is not enough: we need an inverse estimate
that remains uniform as $t\to\infty$. We obtain it by combining the
weighted inverses on the ALE models and the
orbifold. After patching these model solutions, we estimate the error of
the resulting approximate inverse and remove it by a Neumann series. The
result is the uniformly bounded right inverse stated in \cref{lap}.

\subsection{Harmonic analysis for model spaces}

\begin{lemma}\label{EHLaplacian}
	On the ALE K\"{a}hler space $(Q,\omega_Q)$, let $l=1+|z|$ and define the weighted space $C_\delta^{k,\alpha}(Q)$ by using $l$. Then for $-2<\delta<0$, $k\geq0$ and $0<\alpha<1$, the Laplacian with respect to $\omega_Q$, \begin{equation}
		\Delta_Q\colon C_\delta^{k+2,\alpha}(Q)\to C_{\delta-2}^{k,\alpha}(Q)
	\end{equation} is an isomorphism, with a bounded inverse $G_Q$.
\end{lemma}

\begin{proof}
	This is \cite[Theorem 4.5(a)]{joyceAsymptoticallyLocallyEuclidean2001}.
\end{proof}

\begin{lemma}\label{orbiLaplacian}
	On the Calabi--Yau orbifold $(Y,\omega_Y)$, let $(p_i,U_i)$ be the orbifold charts near singularities, define the weighted space $C_\delta^{k,\alpha}(Y^{\mathrm{reg}})$ by using \begin{equation}
		\tilde{r}=\left\{\begin{aligned}
			&r_i,&&\text{in each }U_i \\
			&O(1),&&\text{elsewhere}.
		\end{aligned}\right.
	\end{equation} Then for $-2<\delta<0$, $k\geq0$ and $0<\alpha<1$, for any $f\in C_{\delta-2}^{k,\alpha}(Y^{\mathrm{reg}})$ with \begin{equation}
		\int_Y f\cdot2\Omega_Y\wedge\overline{\Omega}_Y=0,
	\end{equation} there is a unique $u\in C_\delta^{k+2,\alpha}(Y^{\mathrm{reg}})$, such that \begin{equation}
		\Delta_Y u=f,\quad\int_Y u\cdot2\Omega_Y\wedge\overline{\Omega}_Y=0,\quad u=O(r_i^\delta)\text{ near }p_i,
	\end{equation} and that \begin{equation}
		\|u\|_{C_\delta^{k+2,\alpha}(Y^{\mathrm{reg}})}\leq C\|f\|_{C_{\delta-2}^{k,\alpha}(Y^{\mathrm{reg}})}
	\end{equation} with constant $C$ independent of $u,f$. We denote $u=G_Yf$.
\end{lemma}

\begin{proof}
	The proof uses the orbifold Green's function. The Green's kernel $G(x,y)$ of $\Delta_Y$ satisfies \begin{equation}
		|G(x,y)|\leq C(1+\dist_Y(x,y)^{-2}),\quad\int_YG(x,\cdot)\cdot2\Omega_Y\wedge\overline{\Omega}_Y=0.
	\end{equation} Since $\delta-2>-4$, the integral \begin{equation}
		u(x):=\int_Y G(x,\cdot)f\cdot2\Omega_Y\wedge\overline{\Omega}_Y
	\end{equation} is well defined for $x\neq p_i$, with average $0$. Fix $x$ near a singularity and split $U$ into three parts: $r<\frac{r(x)}{2}$, $\frac{r(x)}{2}\leq r\leq2r(x)$, and $2r(x)<r<c_0$. Let $M=\|f\|_{C_{\delta-2}^0(Y^{\mathrm{reg}})}$, then \begin{align}
		&\left|\int_{r<\frac{r(x)}{2}}G(x,\cdot)f\right|\leq Cr(x)^{-2}\int_0^{\frac{r(x)}{2}}Mr^{\delta-2}r^3dr\leq CMr(x)^\delta, \\
		&\left|\int_{\frac{r(x)}{2}\leq r\leq2r(x)}G(x,\cdot)f\right|\leq CMr(x)^{\delta-2}\int_0^{3r(x)}\rho^{-2}\rho^3d\rho\leq CMr(x)^\delta, \\
		&\left|\int_{2r(x)<r<c_0}G(x,\cdot)f\right|\leq C\int_{2r(x)}^{c_0}r^{-2}Mr^{\delta-2}r^3dr\leq CMr(x)^\delta.
	\end{align} The first uses $\delta>-2$ and the last uses $\delta<0$. Away from singularities, $G(x,\cdot)=O(1)$ and $f=O(M)$. Hence in summary, $u(x)=O(r(x)^\delta)$. Then, the H\"{o}lder estimate is from applying interior Schauder estimate for rescaled balls.

	Finally, for uniqueness, if $u_1,u_2$ both satisfy the conditions, then $v=u_1-u_2$ is harmonic, with $v=O(r^\delta)$ near singularities. The elliptic estimate gives $|\nabla v|=O(r^{\delta-1})$. Therefore, for any smooth test function $\phi$, as $r\to0$, \begin{equation}
		\left|\int_{Y\setminus \bigcup_iB(p_i,r)}2v\Delta\phi\right|=\left|-\sum_i\int_{\partial B(p_i,r)}(v\frac{\partial\phi}{\partial\nu}-\phi\frac{\partial v}{\partial\nu})ds\right|\leq C(r^{\delta+2}+r^{\delta+3})\to0.
	\end{equation} So, $v$ is a smooth harmonic function on the orbifold. The zero-average normalization then implies that it vanishes.
\end{proof}

\subsection{Decomposition, patching and parametrix}

For convenience, in this section, we use the metric $\hat{\omega}_a$ and its weighted norm instead of $\omega_t$. Let \begin{equation}
	 w_a:=\sqrt{a}l_{1/a}=\sqrt{a}+\tilde{r}\sim\left\{\begin{aligned}
	 	&\sqrt{a},&&r_i\lesssim\sqrt{a} \\
	 	&r,&&\sqrt{a}\lesssim r_i\lesssim1 \\
	 	&1,&&\tilde{r}\sim1.
	 \end{aligned}\right.
\end{equation} and define $\|\cdot\|_{\hat{C}_{\delta,a}^{k,\alpha}}$ by using $\hat{\omega}_a$ and $w_a$, which relates to the norm in \cref{3.3} by \begin{equation}\label{at}
	\|u\|_{C_{\delta,t}^{k,\alpha}}=t^{-\delta/2}\|u\|_{\hat{C}_{\delta,a}^{k,\alpha}},\quad a=1/t.
\end{equation} And we denote $\hat{C}_{\delta,a}^{k,\alpha,\mathrm{avg}}$ the space \begin{equation}
	\hat{C}_{\delta,a}^{k,\alpha,\mathrm{avg}}=\{f\in\hat{C}_{\delta,a}^{k,\alpha}:\int_X f\hat{\omega}_a^2=0\}.
\end{equation}

We construct the right inverse of $\Delta_{\hat{\omega}_a}\colon\hat{C}_{\delta,a}^{k+2,\alpha}\to\hat{C}_{\delta-2,a}^{k,\alpha,\mathrm{avg}}$. Assume $-2<\delta<0$, for any $f\in\hat{C}_{\delta-2}^{k,\alpha,\mathrm{avg}}$, let $f_i^{\mathrm{in}}=(1-\chi_i)f$ and $f^{\mathrm{out}}=f-\sum_i f_i^{\mathrm{in}}$, where \begin{equation}
	\chi_i=\left\{\begin{aligned}
		&0,&&r_i\leq a^{1/3} \\
		&1,&&r_i\geq 2a^{1/3}
	\end{aligned}\right.
\end{equation} is defined in \cref{3.2}. We choose constants $\beta=\frac{1}{4}$ and $\gamma=\frac{5}{12}$ so that \begin{equation}
	\beta<\frac{1}{3}<\gamma<\frac{1}{2}.
\end{equation} And choose cut-off functions \begin{equation}
	\eta_i^{\mathrm{in}}=\left\{\begin{aligned}
		&0,&&r_i\geq a^\beta \\
		&1,&&r_i\leq 2a^{1/3},
	\end{aligned}\right.\quad\eta^{\mathrm{out}}=\left\{\begin{aligned}
		&0,&&\tilde{r}\leq a^\gamma \\
		&1,&&\tilde{r}\geq a^{1/3}.
	\end{aligned}\right.
\end{equation} For example, we can take \begin{equation}
	\eta_0=\left\{\begin{aligned}
		&0,&&(-\infty,\beta] \\
		&1,&&[\frac{1}{6}+\frac{\beta}{2},\infty)
	\end{aligned}\right.
\end{equation} and let $\eta_i^{\mathrm{in}}=\eta_0(\frac{\log r_i}{\log a})$. A similar choice gives \begin{equation}\label{cutoff}
	r^k|\partial_r^k\eta_i^{\mathrm{in}}|+r^k|\partial_r^k\eta^{\mathrm{out}}|\leq C_k\frac{1}{|\log a|}
\end{equation} for $k\geq 1$.

For the inner part, consider \begin{equation}
	\tau_{a,i}\colon\{|y|<2a^{\beta-1/2}\}\subset Q_i\to\{r_i<2a^\beta\},\quad y\mapsto z=\sqrt{a}y.
\end{equation} Then $\tau_{a,i}$ pulls $\omega_{Q_i,a}$ back to the model metric $a\omega_{Q_i}$ scaled by a factor. Now $\tau_{a,i}^*f_i^{\mathrm{in}}$ is well-defined on $Q_i$, supported on $\{|y|\leq 2a^{-1/6}\}$. We let \begin{equation}
	P_i^{\mathrm{in}}f_i^{\mathrm{in}}:=\eta_i^{\mathrm{in}}(\tau_{a,i}^{-1})^*(aG_{Q_i}(\tau_{a,i}^*f_i^{\mathrm{in}})).
\end{equation}

For the outer part, $f^{\mathrm{out}}$ is viewed as a function on $Y$. Fix a smooth function $k_0$ supported outside the singular charts, with \begin{equation}
	\int_Yk_0\cdot2\Omega_Y\wedge\overline{\Omega}_Y=1.
\end{equation} Then let \begin{equation}
	P^{\mathrm{out}}f^{\mathrm{out}}:=\eta^{\mathrm{out}} G_Y(f^{\mathrm{out}}-k_0\int_Yf^{\mathrm{out}}\cdot 2\Omega_Y\wedge\overline{\Omega}_Y).
\end{equation} And $Pf:=\sum_iP_i^{\mathrm{in}}f_i^{\mathrm{in}}+P^{\mathrm{out}}f^{\mathrm{out}}$.

\begin{proposition}
	For $-2<\delta<0$, \begin{equation}
		\|Pf\|_{\hat{C}_{\delta,a}^{k+2,\alpha}}\leq C\|f\|_{\hat{C}_{\delta-2,a}^{k,\alpha}}.
	\end{equation}
\end{proposition}

\begin{proof}
	Let $M=\|f\|_{\hat{C}_{\delta-2,a}^{k,\alpha}}$. In each $\supp f^{\mathrm{in}}\subset\{r\leq 2a^{1/3}\}$, $w_a\sim\sqrt{a}(1+|y|^4)^{1/4}$. Therefore, by \cref{EHLaplacian}, \begin{equation}
		\|P^{\mathrm{in}}f^{\mathrm{in}}\|_{\hat{C}_{\delta,a}^{k+2,\alpha}}\leq Ca^{-\delta/2}a\|\tau_a^*f^{\mathrm{in}}\|_{C_{\delta-2}^{k,\alpha}(Q)}\leq CM.
	\end{equation}

	In $\supp f^{\mathrm{out}}\subset\{\tilde{r}\geq a^{1/3}\}$, $w_a\sim\tilde{r}$. By the assumption $\int_Xf\hat{\omega}_a^2=0$, we know \begin{equation}
		\begin{aligned}
			\int_Yf^{\mathrm{out}}\cdot2\Omega_Y\wedge\overline{\Omega}_Y&=-\sum_i\int_Xf_i^{\mathrm{in}}\hat{\omega}_a^2+\int_Xf^{\mathrm{out}}(2\Omega\wedge\overline{\Omega}-\hat{\omega}_a^2) \\
			&=-\sum_i\int_{r_i\leq 2a^{1/3}}f_i^{\mathrm{in}}\hat{\omega}_a^2+\sum_i\int_{a^{1/3}\leq r_i\leq 2a^{1/3}}f^{\mathrm{out}}(1-e^F)\cdot2\Omega\wedge\overline{\Omega}+O(Ma^2),
		\end{aligned}
	\end{equation} where $F=\log(\frac{\hat{\omega}_a^2}{2\Omega\wedge\overline{\Omega}})$. Now we estimate \begin{equation}
		\begin{aligned}
			\left|\int_{r\leq 2a^{1/3}}f^{\mathrm{in}}\hat{\omega}_a^2\right|&\leq CM\int_{r\leq 2a^{1/3}}w_a^{\delta-2}\hat{\omega}_a^2 \\
			&\leq CM(\int_0^{\sqrt{a}}a^{\frac{\delta-2}{2}}r^3dr+\int_{\sqrt{a}}^{2a^{1/3}}r^{\delta-2}r^3dr) \\
			&\leq CMa^{\frac{\delta+2}{3}}.
		\end{aligned}
	\end{equation} And that \begin{equation}
		\left|\int_{a^{1/3}\leq r\leq 2a^{1/3}}f^{\mathrm{out}}(1-e^F)\cdot2\Omega\wedge\overline{\Omega}\right|\leq CM\int_{a^{1/3}}^{2a^{1/3}}r^{\delta-2}a^{2/3}r^3dr\leq CMa^{\frac{\delta+4}{3}}
	\end{equation} by using \cref{F}. Since $\delta>-2$, all $a^{\frac{\delta+2}{3}}$, $a^{\frac{\delta+4}{3}}$ and $a^2$ are small as $a\ll1$. So we find \begin{equation}
		\|P^{\mathrm{out}}f^{\mathrm{out}}\|_{\hat{C}_{\delta,a}^{k+2,\alpha}}\leq C\left\|f^{\mathrm{out}}-k_0\int_Y f^{\mathrm{out}}\cdot2\Omega_Y\wedge\overline{\Omega}_Y\right\|_{C_{\delta-2}^{k,\alpha}(Y^{\mathrm{reg}})}\leq CM.
	\end{equation} In summary, $\|Pf\|_{\hat{C}_{\delta,a}^{k+2,\alpha}}\leq CM$.
\end{proof}

\begin{proposition}
	For $-2<\delta<0$, $\Delta=\Delta_{\hat{\omega}_a}$, \begin{equation}
		\|\Delta Pf-f\|_{\hat{C}_{\delta-2,a}^{k,\alpha}}\leq C\frac{1}{|\log a|}\|f\|_{\hat{C}_{\delta-2,a}^{k,\alpha}}.
	\end{equation}
\end{proposition}

\begin{proof}
	Let $M=\|f\|_{\hat{C}_{\delta-2,a}^{k,\alpha}}$, and let $\Delta_{Q,a}$ be the inner Laplacian, $\Delta_{Q,a}((\tau_a^{-1})^*(aG_Q(\tau_a^*f^{\mathrm{in}})))=f^{\mathrm{in}}$. Then \begin{equation}
		\begin{aligned}
			\Delta Pf-f&=\sum_i(\Delta P_i^{\mathrm{in}}f_i^{\mathrm{in}}-f_i^{\mathrm{in}})+(\Delta P^{\mathrm{out}}f^{\mathrm{out}}-f^{\mathrm{out}}) \\
			&=\sum_i\eta_i^{\mathrm{in}}(\Delta-\Delta_{Q_i,a})\left((\tau_{a,i}^{-1})^*(aG_{Q_i}(\tau_{a,i}^*f_i^{\mathrm{in}}))\right) \\
			&\phantom{=}+\eta^{\mathrm{out}}(\Delta-\Delta_Y)\left(G_Y(f^{\mathrm{out}}-k_0\int_Yf^{\mathrm{out}}\cdot2\Omega_Y\wedge\overline{\Omega}_Y)\right) \\
			&\phantom{=}-k_0\int_Yf^{\mathrm{out}}\cdot2\Omega_Y\wedge\overline{\Omega}_Y \\
			&\phantom{=}+\sum_i\left(\Delta\eta_i^{\mathrm{in}}\cdot(\tau_{a,i}^{-1})^*(aG_{Q_i}(\tau_{a,i}^*f_i^{\mathrm{in}}))+\langle\nabla\eta_i^{\mathrm{in}},\nabla(\tau_{a,i}^{-1})^*(aG_{Q_i}(\tau_{a,i}^*f_i^{\mathrm{in}}))\rangle\right) \\
			&\phantom{=}+\left(\Delta\eta^{\mathrm{out}}\cdot G_Y(f^{\mathrm{out}}-k_0\int_Yf^{\mathrm{out}}\cdot2\Omega_Y\wedge\overline{\Omega}_Y)+\langle\nabla\eta^{\mathrm{out}},\nabla G_Y(f^{\mathrm{out}}-k_0\int_Yf^{\mathrm{out}}\cdot2\Omega_Y\wedge\overline{\Omega}_Y)\rangle\right) \\
			&=A_1+A_2-k_0\int_Yf^{\mathrm{out}}\cdot2\Omega_Y\wedge\overline{\Omega}_Y+A_3+A_4.
		\end{aligned}
	\end{equation}

	$A_1$ is supported on $\{a^{1/3}\leq r\leq a^\beta\}$. Here, the metric difference \begin{equation}
		\omega_{Q,a}-\omega_{\mathrm{Euc}}=\sqrt{-1}\partial\bar{\partial}(a\delta_{-\frac{1}{2}\log a}^*\Psi_Q-r^2),
	\end{equation} so \begin{equation}
		r^k|\nabla^k(\omega_{Q,a}-\omega_{\mathrm{Euc}})|_{\omega_{\mathrm{Euc}}}\leq C\frac{a^2}{r^4}\leq Ca^{2/3},
	\end{equation} and $\omega_{Q,a}$ is comparable to $\omega_{\mathrm{Euc}}$. Again, \begin{equation}
		\hat{\omega}_a-\omega_{Q,a}=\sqrt{-1}\partial\bar{\partial}(\chi(u_{Y,a}-(a\delta_{-\frac{1}{2}\log a}^*\Psi_Q-r^2))),
	\end{equation} we get \begin{equation}
		r^k|\nabla^k(\hat{\omega}_a-\omega_{Q,a})|_{\omega_{Q,a}}\leq C(r^2+\frac{a^2}{r^4})\leq Ca^{2\beta}.
	\end{equation} Therefore, from the metric discrepancy \begin{equation}
		\|A_1\|_{\hat{C}_{\delta-2,a}^{k,\alpha}}\leq Ca^{2\beta}\|(\tau_a^{-1})^*(aG_Q(\tau_a^*f_1))\|_{\hat{C}_{\delta,a}^{k+2,\alpha}}\leq CMa^{2\beta}.
	\end{equation} The same argument shows \begin{equation}
		\|A_2\|_{\hat{C}_{\delta-2,a}^{k,\alpha}}\leq CM(a^{2-4\gamma}+a^{2/3}+a).
	\end{equation}

	From \labelcref{cutoff}, it's easy to obtain \begin{equation}
		\|A_3\|_{\hat{C}_{\delta-2,a}^{k,\alpha}}+\|A_4\|_{\hat{C}_{\delta-2,a}^{k,\alpha}}\leq CM\frac{1}{|\log a|}.
	\end{equation} Finally, recall that \begin{equation}
		\left|\int_Yf^{\mathrm{out}}\cdot2\Omega_Y\wedge\overline{\Omega}_Y\right|\leq CMa^{\frac{\delta+2}{3}},
	\end{equation} thus the conclusion.
\end{proof}

So we can choose $a\ll 1$ even smaller, such that $\|\Delta P-1\|\leq\frac{C}{|\log a|}<1$. This implies \begin{equation}
	R:=P\sum_{i=0}^\infty(1-\Delta P)^i
\end{equation} is well-defined on $\hat{C}_{\delta-2,a}^{k,\alpha,\mathrm{avg}}$, and is a right inverse of $\Delta$, with uniformly bounded operator norm \begin{equation}
	\|R\|\leq\|P\|\sum_{i=0}^\infty\|1-\Delta P\|^i\leq C.
\end{equation} Back to $\omega_t$, use \labelcref{at} and $\Delta_{\omega_t}=t^{-1}\Delta$, we find $u=G_tf:=R(tf)$ solves $\Delta_{\omega_t}u=f$, with bounded norm \begin{equation}
	\|u\|_{C_{\delta,t}^{k+2,\alpha}}\leq C\|f\|_{C_{\delta-2,t}^{k,\alpha}}.
\end{equation} We've proved the following conclusion, where \begin{equation}
	C_{\delta,t}^{k,\alpha,\mathrm{avg}}=\{f\in C_{\delta,t}^{k,\alpha}:\int_X f\omega_t^2=0\}.
\end{equation}

\begin{proposition}\label{lap}
	For $-2<\delta<0$, $k\geq0$ and $0<\alpha<1$, there exists a right inverse $G_t$ of \begin{equation}
		\Delta_{\omega_t}\colon C_{\delta,t}^{k+2,\alpha}\to C_{\delta-2,t}^{k,\alpha,\mathrm{avg}},
	\end{equation} whose operator norm $\|G_t\|$ is uniformly bounded as $t\gg1$.
\end{proposition}

\section{Parameter dependence and iteration}\label{s5}

The right inverse from \cref{s4} permits corrections on each slice, but a
correction of the full ansatz must also be controlled in the parameter $t$.
In particular, changing $\omega_t$ by
$\sqrt{-1}\partial\bar\partial u_t$ requires changing $h_t$ by
$-\frac12\partial_t^2u_t$. We therefore establish stability of the inverse
estimates under small perturbations and prove estimates for the parameter
derivatives of normalized Poisson solutions. These results justify a finite
correction scheme that improves the Calabi--Yau error while preserving the
reduced cohomology class and positivity. Our goal is sufficient decay for applying
the Tian--Yau--Hein package, rather than an exact solution of the reduced
equations.

In this section, we assume $(\omega_t,h_t)$ is a general family of K\"{a}hler metrics that satisfies the approximate Gibbons--Hawking ansatz \labelcref{agh}, and denote $(\omega_t^{(0)},h_t^{(0)})$ the gluing metric constructed in \cref{3.2}.

\begin{lemma}\label{perlap}
	Let $-2<\delta<0$, assume $\Delta_{\omega_t}\colon C_{\delta,t}^{k+2,\alpha}\to C_{\delta-2,t}^{k,\alpha,\mathrm{avg}}$ has a uniformly bounded right inverse $G_t$, and $\omega_t^2\leq C\cdot2t^2\Omega\wedge\overline{\Omega}$ uniformly. If we perturb $\omega_t'=\omega_t+\sqrt{-1}\partial\bar{\partial}u_t$ by $u_t$ with small error in the sense \begin{equation}
		\|u_t\|_{C_{\delta,t}^{k+2,\alpha}}\leq\epsilon(t)\to0.
	\end{equation} Denote $C_{\delta-2,t}^{k,\alpha,\mathrm{avg}'}$ the zero-average space with respect to $(\omega_t')^2$. Then for sufficiently large $t$, the Laplacian $\Delta_{\omega_t'}\colon C_{\delta,t}^{k+2,\alpha}\to C_{\delta-2,t}^{k,\alpha,\mathrm{avg}'}$ also has a uniformly bounded right inverse $G_t'$.
\end{lemma}

\begin{proof}
	Let \begin{equation}
		C_{\delta,t,\mathrm{avg}}^{k+2,\alpha}=\{u\in C_{\delta,t}^{k+2,\alpha}:\int_Xu\cdot2\Omega\wedge\overline{\Omega}=0\}.
	\end{equation} be the zero-average space with respect to $2\Omega\wedge\overline{\Omega}$. We denote the projection \begin{equation}
		Qf=f-\frac{\int_X f\cdot2\Omega\wedge\overline{\Omega}}{\int_X 2\Omega\wedge\overline{\Omega}},
	\end{equation} and let that $\widetilde{\Delta}_t=Q\Delta_{\omega_t}$, $\widetilde{\Delta}_t'=Q\Delta_{\omega_t'}$. Also we let \begin{equation}
		\widetilde{G}_tf=QG_t(f-\frac{\int_Xf\omega_t^2}{\int_X\omega_t^2}).
	\end{equation} For $\sigma=\delta$ or $\delta-2$, since \begin{equation}
		|Qf|\leq |f|+\frac{\int_X l_t^\sigma\cdot2\Omega\wedge\overline{\Omega}}{\int_X2\Omega\wedge\overline{\Omega}}\cdot\sup_X(l_t^{-\sigma}|f|),
	\end{equation} where \begin{equation}
		\int_X l_t^\sigma\cdot2\Omega\wedge\overline{\Omega}\leq C(\int_0^{t^{-1/2}}r^3dr+\int_{t^{-1/2}}^1t^{\sigma/2}r^\sigma r^3dr+t^{\sigma/2})\leq Ct^{\sigma/2}\leq Cl_t^\sigma.
	\end{equation} We know that $Q$ is uniformly bounded on $C_{\sigma,t}^{k,\alpha}$. Similarly, the projection \begin{equation}
		f\mapsto f-\frac{\int_Xf\omega_t^2}{\int_X\omega_t^2}
	\end{equation} is uniformly bounded on $C_{\delta-2,t}^{k,\alpha}$ because \begin{equation}
		\int_X l_t^{\delta-2}\omega_t^2\leq Ct^2\int_Xl_t^{\delta-2}\cdot2\Omega\wedge\overline{\Omega},\quad\int_X\omega_t^2\sim t^2\int_X2\Omega\wedge\overline{\Omega}.
	\end{equation} Therefore, the operators $\widetilde{\Delta}_t$, $\widetilde{\Delta}_t'$ and $\widetilde{G}_t$ are uniformly bounded.

	Now we restrict to the space $C_{\delta-2,t,\mathrm{avg}}^{k,\alpha}$, then $\widetilde{\Delta}_t\widetilde{G}_t=1_{C_{\delta-2,t,\mathrm{avg}}^{k,\alpha}}$ is the identity map. Therefore, from the equality \begin{equation}
		\widetilde{\Delta}_t'=\widetilde{\Delta}_t(1-\widetilde{G}_t(\widetilde{\Delta}_t-\widetilde{\Delta}_t')),
	\end{equation} formally we have \begin{equation}
		\widetilde{G}_t'=(1-\widetilde{G}_t(\widetilde{\Delta}_t-\widetilde{\Delta}_t'))^{-1}\widetilde{G}_t.
	\end{equation} The condition of $u_t$ controls the metric discrepancy, which implies that \begin{equation}
		\|(\widetilde{\Delta}_t-\widetilde{\Delta}_t')v\|_{C_{\delta-2,t}^{k,\alpha}}\leq C\epsilon(t)\|v\|_{C_{\delta,t}^{k+2,\alpha}},
	\end{equation} so the operator norm \begin{equation}
		\|\widetilde{G}_t(\widetilde{\Delta}_t-\widetilde{\Delta}_t')\|_{C_{\delta,t,\mathrm{avg}}^{k+2,\alpha}\to C_{\delta,t,\mathrm{avg}}^{k+2,\alpha}}\leq C\epsilon(t)\to0.
	\end{equation} Provided that $t$ is sufficiently large, the inverse operator \begin{equation}
		(1-\widetilde{G}_t(\widetilde{\Delta}_t-\widetilde{\Delta}_t'))^{-1}=\sum_{i=0}^\infty(\widetilde{G}_t(\widetilde{\Delta}_t-\widetilde{\Delta}_t'))^i,
	\end{equation} and thus $\widetilde{G}_t'$, is well-defined. Finally, $G_t'f=\widetilde{G}_t'Qf$ for function $f$ with $\int_Xf(\omega_t')^2=0$.
\end{proof}

Now we study the parameter dependence of the solution. We say $\omega_t$ is uniformly comparable to $\omega_t^{(0)}$, if $\frac{1}{C}\omega_t^{(0)}\leq\omega_t\leq C\omega_t^{(0)}$ for a constant $C$ independent of $t$. For example, if $\omega_t=\omega_t^{(0)}+\sqrt{-1}\partial\bar{\partial}u_t$ with $u_t$ small in the weight norms.

\begin{definition}
	We say $\omega_t$ satisfies the \textbf{coefficients condition}, if $\omega_t$ is uniformly comparable to $\omega_t^{(0)}$, and \begin{equation}
		l_t^it^j|\nabla_{\omega_t^{(0)}}^i\partial_t^j\omega_t|_{\omega_t^{(0)}}\leq C_{i,j}
	\end{equation} for every $i,j\geq0$.
\end{definition}

\begin{lemma}\label{tdep}
	Let $-2<\delta<0$, assume $\omega_t$ satisfies the coefficients condition, and $\Delta_{\omega_t}\colon C_{\delta,t}^{k+2,\alpha}\to C_{\delta-2,t}^{k,\alpha,\mathrm{avg}}$ admits a uniformly bounded right inverse $G_t$, modified by $\widetilde{G}_t=QG_t$ for having average $0$. Then for every $m\geq0$, \begin{equation}
		\|\widetilde{G}_t f\|_{C_{\delta,t}^{k+2,\alpha;m}}\leq C\|f\|_{C_{\delta-2,t}^{k,\alpha,m}}
	\end{equation} with $C$ independent of $t$.
\end{lemma}

\begin{proof}
	Suppose $\Delta_{\omega_t}u_t=f_t$. Take derivatives in $t$, we have \begin{equation}
		\Delta_t(t^m\partial_t^m u_t)=t^m\partial_t^m f_t-\sum_{i=1}^m\binom{m}{i}(t^i\partial_t^i\Delta_{\omega_t})(t^{m-i}\partial_t^{m-i}u_t).
	\end{equation} The zero-average condition, $\int_Xu_t\cdot2\Omega\wedge\overline{\Omega}=0$, holds when taking derivatives in $t$. Thus, \begin{equation}
		t^m\partial_t^mu_t=\widetilde{G}_t(t^m\partial_t^m f_t-\sum_{i=1}^m\binom{m}{i}(t^i\partial_t^i\Delta_{\omega_t})(t^{m-i}\partial_t^{m-i}u_t)).
	\end{equation} Now the coefficients condition controls the coefficients of $t^i\partial_t^i\Delta_{\omega_t}$, so that \begin{equation}
		\|(t^i\partial_t^i\Delta_{\omega_t})v\|_{C_{\delta-2,t}^{k,\alpha}}\leq C\|v\|_{C_{\delta,t}^{k+2,\alpha}}.
	\end{equation} We get \begin{equation}
		t^m\|\partial_t^mu_t\|_{C_{\delta,t}^{k+2,\alpha}}\leq C(t^m\|\partial_t^mf_t\|_{C_{\delta-2,t}^{k,\alpha}}+\sum_{i=1}^mt^{m-i}\|\partial_t^{m-i}u_t\|_{C_{\delta,t}^{k+2,\alpha}}),
	\end{equation} and the desired estimate is proved by induction.
\end{proof}

The parameter dependence holds in perturbation, provided that the perturbation is small in the sense of $\|\cdot\|_{C_{\delta,t}^{k+2,\alpha;m}}$.

\begin{lemma}\label{pertdep}
	Let $-2<\delta<0$, assume $\omega_t$ satisfies the coefficients condition, $\omega_t'=\omega_t+\sqrt{-1}\partial\bar{\partial}u_t$, with $u_t$ small in the sense \begin{equation}
		\|u_t\|_{C_{\delta,t}^{k+2,\alpha;m}}\leq C_{k,\alpha,m}\epsilon(t)\to0
	\end{equation} for every $k\geq0$, $0<\alpha<1$, $m\geq0$. Then for sufficiently large $t$, $\omega_t'$ also satisfies the coefficients condition.
\end{lemma}

\begin{proof}
	The condition of $u_t$ implies \begin{equation}\label{6}
		|\nabla_{\omega_t^{(0)}}^i\partial_t^j(\omega_t-\omega_t')|_{\omega_t^{(0)}}\leq C\epsilon(t)l_t^{-i-2+\delta}t^{-j}.
	\end{equation} Since $\delta-2<0$, $l_t\geq1$, and $\epsilon(t)\to0$, $\epsilon(t)l_t^{-2+\delta}$ is small and thus \begin{equation}\label{7}
		|\nabla_{\omega_t^{(0)}}^i\partial_t^j\omega_t'|_{\omega_t^{(0)}}\leq Cl_t^{-i}t^{-j}.
	\end{equation}
\end{proof}

We now use these estimates to perform one correction step. For a family $(\omega_t,h_t)$ satisfying \labelcref{agh}, let the error term \begin{equation}
	E_t=\frac{2h_t\,2\Omega\wedge\overline{\Omega}}{\omega_t^2}-1,
\end{equation}
the integral normalization in \labelcref{agh} gives
$\int_X E_t\,\omega_t^2=0$. Thus $E_t$ satisfies the compatibility condition
for the Poisson equation. We fix the additive constant in its solution
using the $t$-independent measure $2\Omega\wedge\overline{\Omega}$; this
allows the normalization to be differentiated without introducing
derivatives of a varying volume form. For sufficiently large $t$, we solve $u_t$ from the equation \begin{equation}\label{9}
	\left\{\begin{aligned}
		&\Delta_{\omega_t}u_t=E_t \\
		&\int_X u_t\cdot2\Omega\wedge\overline{\Omega}=0.
	\end{aligned}\right.
\end{equation} Then, \begin{equation}\label{10}
	\omega_t'=\omega_t+\sqrt{-1}\partial\bar{\partial}u_t,\quad h_t'=h_t-\frac{1}{2}\partial_t^2u_t
\end{equation} defines a new family satisfying \labelcref{agh}, except possibly for positivity. The new error term \begin{equation}\label{newerror}
	E_t'=\frac{2h_t'\cdot2\Omega\wedge\overline{\Omega}}{(\omega_t')^2}-1=-\frac{\frac{\partial_t^2u_t}{2h_t}(1+E_t)+\frac{(\sqrt{-1}\partial\bar{\partial }u_t)^2}{\omega_t^2}}{1+E_t+\frac{(\sqrt{-1}\partial\bar{\partial}u_t)^2}{\omega_t^2}}.
\end{equation}

\begin{proposition}\label{it}
	Let $-2<\delta<0$, and assume \begin{itemize}
		\item $\omega_t$ satisfies the coefficients condition;
		\item $\Delta_{\omega_t}\colon C_{\delta,t}^{k+2,\alpha}\to C_{\delta-2,t}^{k,\alpha,\mathrm{avg}}$ has a uniformly bounded right inverse $G_t$;
		\item $\|2h_t/t^2\|_{C_{0,t}^{k,\alpha;m}}\leq C_{k,\alpha,m}$;
		\item $\|E_t\|_{C_{\delta-2,t}^{k,\alpha;m}}\leq C_{k,\alpha,m}\epsilon(t)\to0$;
	\end{itemize} hold for every $k\geq0$, $0<\alpha<1$, $m\geq0$. Then for sufficiently large $t\gg1$, the ansatz $(\omega_t',h_t')$ in \labelcref{9,10} is well-defined, $\omega_t'$ satisfies the coefficients condition, and that \begin{equation}
		\|\frac{2h_t'}{t^2}\|_{C_{0,t}^{k,\alpha;m}}\leq A_{k,\alpha,m},\quad\|E_t'\|_{C_{\delta-2,t}^{k,\alpha;m}}\leq B_{k,\alpha,m}(\epsilon(t)^2+t^{-3}\epsilon(t)),
	\end{equation} for constants $A_{k,\alpha,m},B_{k,\alpha,m}$ uniform in $t$.
\end{proposition}

\begin{proof}
	By \cref{tdep} we know $\|u_t\|_{C_{\delta,t}^{k+2,\alpha;m+2}}\leq C\|E_t\|_{C_{\delta-2,t}^{k,\alpha;m+2}}\leq C\epsilon(t)\to0$. In particular, $\omega_t'$ is K\"{a}hler for $t\gg1$ and satisfies the coefficients condition. Meanwhile, the coefficients condition of $\omega_t$ controls $\|\frac{\omega_t^2}{(\omega_t^{(0)})^2}\|_{C_{0,t}^{k,\alpha;m}}$. These give \begin{equation}\label{11}
		\left\|\frac{(\sqrt{-1}\partial\bar{\partial}u_t)^2}{\omega_t^2}\right\|_{C_{\delta-2,t}^{k,\alpha;m}}\leq C\left\|\frac{(\sqrt{-1}\partial\bar{\partial}u_t)^2}{(\omega_t^{(0)})^2}\right\|_{C_{\delta-2,t}^{k,\alpha;m}}\leq C\sup_Xl_t^{\delta-2}\|u_t\|_{C_{\delta,t}^{k+2,\alpha;m}}^2\leq C\epsilon(t)^2.
	\end{equation} Since \begin{equation}
		t^k\partial_t^k(\frac{\partial_t^2u_t}{2h_t})=\sum_{i=0}^k\binom{k}{i}(t^i\partial_t^i(2h_t)^{-1})(t^{k-i}\partial_t^{k-i+2}u_t),
	\end{equation} we find \begin{equation}\label{12}
		\|\frac{\partial_t^2u_t}{2h_t}\|_{C_{\delta-2,t}^{k,\alpha;m}}\leq Ct^{-2}\cdot t^{-2}\sup_Xl_t^2\|u_t\|_{C_{\delta,t}^{k,\alpha;m+2}}\leq Ct^{-3}\epsilon(t)\to0.
	\end{equation} Here we use $l_t\leq C\sqrt{t}$ on $X$. This ensures $h_t'=h_t-\frac{1}{2}\partial_t^2u_t>0$ for $t\gg1$. Together with \labelcref{newerror,11,12}, it's easy to show \begin{equation}
		\|E_t'\|_{C_{\delta-2,t}^{k,\alpha;m}}\leq C'(\epsilon(t)^2+t^{-3}\epsilon(t)).
	\end{equation} And \begin{equation}
		\|\frac{2h_t'}{t^2}\|_{C_{0,t}^{k,\alpha;m}}=\|\frac{2h_t}{t^2}(1-\frac{\partial_t^2u_t}{2h_t})\|_{C_{0,t}^{k,\alpha;m}}\leq C.
	\end{equation}
\end{proof}

Now we can apply the iteration process. Start from $\omega_t^{(0)}$, and let $F_t^{(0)}=\log(\frac{(\omega_t^{(0)})^2}{t^2\cdot2\Omega\wedge\overline{\Omega}})$. Then initial error term \begin{equation}
	E_t^{(0)}=\frac{2h_t^{(0)}\cdot2\Omega\wedge\overline{\Omega}}{(\omega_t^{(0)})^2}-1=(\frac{2h_t^{(0)}}{t^2}-1)e^{-F_t^{(0)}}+(e^{-F_t^{(0)}}-1).
\end{equation} In $\{t^{-1/3}\leq r\leq2t^{-1/3}\}$, where $l_t\sim t^{1/6}$, so \cref{F} implies \begin{equation}
	\|F_t^{(0)}\|_{C_{\delta-2,t}^{k,\alpha;m}(\{t^{-1/3}\leq r\leq2t^{-1/3}\})}\leq Ct^{-2/3}t^{-\frac{\delta-2}{6}}\leq Ct^{-\frac{2+\delta}{6}}.
\end{equation} In $\{\tilde{r}>2t^{-1/3}\}$, since $l_t\leq C\sqrt{t}$, \begin{equation}
	F_t^{(0)}=\log(1+\frac{[\gamma]^2}{4\pi^2c_1(L)^2}\frac{1}{t^2}),\quad\|F_t^{(0)}\|_{C_{\delta-2,t}^{k,\alpha;m}(\{\tilde{r}>2t^{-1/3}\})}\leq Ct^{-2}t^{-\frac{\delta-2}{2}}\leq Ct^{-\frac{2+\delta}{2}}.
\end{equation} Use \cref{h}, there is \begin{equation}
	\|\frac{2h_t^{(0)}}{t^2}-1\|_{C_{\delta-2,t}^{k,\alpha;m}}\leq Ct^{-2}t^{-\frac{\delta-2}{2}}\leq Ct^{-\frac{2+\delta}{2}}.
\end{equation} For $-2<\delta<0$, these estimates give \begin{equation}
	\|E_t^{(0)}\|_{C_{\delta-2,t}^{k,\alpha;m}}\leq Ct^{-\frac{2+\delta}{6}}.
\end{equation} Therefore, \cref{coe,lap,F,h} guarantee the conditions of \cref{it}. So we solve $u_t^{(0)}$ from the error term as in \labelcref{9}, and modify $\omega_t^{(1)}=\omega_t^{(0)}+\sqrt{-1}\partial\bar{\partial}u_t^{(0)}$. \Cref{perlap,tdep,pertdep,it} ensure this process can be made any prescribed finite number of times, as we enlarge $t\gg1$ suitably. So, start from $\epsilon_0(t)=t^{-a_0}$, the iteration improves the errors by \begin{equation}
	\epsilon_n(t)=t^{-a_n},\quad a_n=\min(2a_{n-1},a_{n-1}+3).
\end{equation} We take \begin{equation}
	\delta=-\frac{1}{8},\quad a_0=\frac{2+\delta}{6}=\frac{5}{16},
\end{equation} after four iterations, the error \begin{equation}
	a_4=5,\quad\epsilon_4(t)=t^{-5}.
\end{equation}

The exponent $5$ is chosen to meet the decay requirement in the final
existence theorem. \Cref{r} will show that the distance from a fixed
base point is comparable to $t^2$. Consequently, the $O(t^{-5})$ error
obtained after four correction steps gives the $O(r^{-5/2})$
Ricci potential decay used in \cref{s7}. We record the improved family
and its weighted estimates in the following theorem.

\begin{theorem}\label{metric}
	Let $\delta=-1/8$. There exists $u_t$, satisfying $\|u_t\|_{C_{\delta,t}^{k+2,\alpha;m}}\leq Ct^{-5/16}$ for every $k\geq0$, $0<\alpha<1$, $m\geq0$, and $\omega_t=\omega_t^{(0)}+\sqrt{-1}\partial\bar{\partial}u_t$, $h_t=h_t^{(0)}-\frac{1}{2}\partial_t^2u_t$, such that \begin{equation}
		\left\|\frac{2h_t\cdot2\Omega\wedge\overline{\Omega}}{\omega_t^2}-1\right\|_{C_{\delta-2,t}^{k,\alpha;m}}\leq Ct^{-5}.
	\end{equation} In particular, we have \begin{equation}
		\left|\frac{2h_t\cdot2\Omega\wedge\overline{\Omega}}{\omega_t^2}-1\right|\leq Ct^{-5}l_t^{\delta-2}\leq Ct^{-5}.
	\end{equation}
\end{theorem}

\section{Transfer to divisor complement}\label{s6}

The improved family from \cref{s5} determines a model metric on the normal
bundle. We now turn this model into a global K\"{a}hler metric on $M$ in the
prescribed class. This requires comparing the complex structures under a
tubular identification, gluing the transferred end metric to an interior
K\"{a}hler metric, and imposing the integral compatibility condition for the
final Monge--Amp\`{e}re equation. The global potential from \cref{s2} allows
the gluing to be performed within the prescribed class. The remaining
additive constant in the relation between $s$ and $t$ will then be used to
satisfy the integral condition.

\subsection{The global gluing}

We now return to the setting in \cref{intro}. Let $X=D$, we'll construct a smooth K\"{a}hler metric on $M=\overline{M}\setminus D$.

We apply the method in \cref{complement}. And use $\alpha_0,\beta_0$ defined in \labelcref{13,18}, \begin{equation}
	\omega_t=\omega_t^{(0)}+\sqrt{-1}\partial\bar{\partial}u_t=t\alpha_0+\beta_0+\sqrt{-1}\partial\bar{\partial}(v_t+u_t).
\end{equation} Since the corresponding cohomology classes agree, set \begin{equation}
	\alpha_0=\alpha|_D+\sqrt{-1}\partial\bar{\partial}\phi_\alpha,\quad\beta_0=\beta|_D+\sqrt{-1}\partial\bar{\partial}\phi_\beta,
\end{equation} we know $\omega_t=t\alpha|_D+\beta|_D+\sqrt{-1}\partial\bar{\partial}\phi_t$ for \begin{equation}
	\phi_t=t\phi_\alpha+\phi_\beta+v_t+u_t.
\end{equation} Normalize it to have average $0$, \begin{equation}
	\tilde{\phi}_t=\phi_t-\frac{\int_X \phi_t\cdot2\Omega\wedge\overline{\Omega}}{\int_X 2\Omega\wedge\overline{\Omega}}.
\end{equation} Here $\phi_\alpha,\phi_\beta$ are fixed functions independent of $t$, together with \cref{phi,metric}, these estimates show \begin{equation}\label{tildephi}
	t^k\|\partial_t^k\tilde{\phi}_t\|_{C^0}\leq C_kt
\end{equation} for $k\geq0$. In particular, $\|\partial_t\tilde{\phi}_t\|_{C^0}\leq C$. So, putting $n=2$ and $[\alpha|_D],[\beta|_D]$ into \labelcref{s} gives the relation \begin{equation}\label{st}
	s=\frac{t^3}{6}+\frac{[\beta|_D]^2}{8\pi^2c_1(L)^2}t-\frac{1}{2}\partial_t\tilde{\phi}_t+c=\frac{t^3}{6}+\frac{[\beta|_D]^2}{8\pi^2c_1(L)^2}t+O(1).
\end{equation} In particular, the metric $\omega_L$ defined by $(\omega_t,h_t)$ and \cref{potential} is defined on a punctured tubular neighborhood $N=\{t\gg1\}\subset L^\times$ of the zero section. We take it as the model metric.

Choose an ambient Hermitian metric on $\overline{M}$, then the normal exponential map gives a diffeomorphism $\Phi$ from $N$ to a puncture tubular neighborhood of $D$ in $\overline{M}$. This diffeomorphism need not be holomorphic, so merely transporting
$\omega_L$ does not produce the desired K\"{a}hler form on $M$. Instead, we
transfer its potential and apply the complex differential operators on
$M$. The following comparison estimates show that this changes the model
only by an exponentially small error, which is negligible relative to the
$O(t^{-5})$ error obtained in \cref{s5}.

\begin{proposition}\label{cstr}
	Let $J_{\overline{M}}$ and $J_L$ be the complex structures on $\overline{M}$ and $L$. Let $\Omega_M$ and $\Omega_L$ be the holomorphic volume forms on $M$ and $L$, related to the holomorphic volume form $\Omega_D$ on $D$ by \labelcref{OmegaL} and $\operatorname{Res}_D \Omega_M=\Omega_D$. Let $\pi_L\colon L\to D$ and $\beta$ be defined as above. Then \begin{equation}
		|\nabla_{g_L}^k(\Phi^*J_{\overline{M}}-J_L)|_{g_L}+|\nabla_{g_L}^k(\Phi^*\Omega_M-\Omega_L)|_{g_L}+|\nabla_{g_L}^k(\Phi^*\beta-\pi_L^*(\beta|_D))|_{g_L}=O(e^{-(1-\epsilon)s})=O(e^{-ct^3})
	\end{equation} for all $k\geq0$, $\epsilon>0$.
\end{proposition}

\begin{proof}
	The proof is essentially in \cite[Proposition 3.4]{heinNilpotentStructuresCollapsing2021} and \cite[Proposition 4.3]{chenCalabiYauMetricsCalabi2024}. Note that their proof doesn't require the ampleness of $L$. Let $(z_1,z_2,w)$ be local holomorphic coordinates on $\overline{M}$ such that $D$ is cut out by  $w=0$, also viewed as holomorphic coordinates on $L$. And we follow the notations in \cref{complement}. Instead of the Calabi model space, our metric \begin{equation}
		g_L=h_t^{-1}\Theta^2+h_tdt^2+g_t.
	\end{equation}

	We fix a large $T$, by integrating from $|\partial_tg_t^{(0)}|_{g_t^{(0)}}\leq Ct^{-1}$, we know \begin{equation}
		(\frac{t}{T})^{-C}g_T^{(0)}\leq g_t^{(0)}\leq(\frac{t}{T})^Cg_T^{(0)}.
	\end{equation} So the coefficients of $g_t^{(0)}$ with respect to the fixed $g_T^{(0)}$, as well as the inverse, are bounded by a polynomial in $t$. Let $S_t=\nabla_{g_t^{(0)}}-\nabla_{g_T^{(0)}}$, then \begin{equation}
		\partial_tS_t=\frac{1}{2}(g_t^{(0)})^{-1}\ast(\nabla_{g_t^{(0)}}\partial_t g_t^{(0)}).
	\end{equation} Expand $\nabla_{g_T^{(0)}}^k\partial_tS_t=(\nabla_{g_t^{(0)}}-S_t)^k\partial_tS_t$, integrate in $t$, and then induce on the order $k$, this gives the polynomial bound for every derivative of the coefficients of $g_t^{(0)}$. Then the coefficients condition gives the polynomial bound for $g_t$'s coefficients, inverse, and derivatives. Along the fiber directions, $h_t$ also gives a polynomial bound. That is, \begin{equation}\label{15}
		|\nabla_{g_L}^k z_i|_{g_L}\leq C_kt^{N_k},\quad|\nabla_{g_L}^kt|_{g_L}\leq C_kt^{N_k}.
	\end{equation}

	Now \begin{equation}
		\log w=-s+\frac{\rho}{2}+\sqrt{-1}\sigma=-\frac{t^3}{6}-\frac{[\beta|_D]^2}{8\pi^2c_1(L)^2}t+\frac{1}{2}\partial_t\tilde{\phi}_t-c+\frac{\rho}{2}+\sqrt{-1}\sigma.
	\end{equation} The radial derivatives of $\partial_t\tilde{\phi}_t$ is controled by \cref{tildephi}. For mixed derivatives, recall \begin{equation}
		\partial_t^2\tilde{\phi}_t=t^2+\frac{[\beta|_D]^2}{4\pi^2c_1(L)^2}-2h_t.
	\end{equation} Apply $\nabla_{g_T^{(0)}}^k$ and integrate on $t$, we know the derivatives of $\partial_t\tilde{\phi}_t$ develops a polynomial bound as well. As a result, \begin{equation}
		|\nabla_{g_L}^k\log w|_{g_L}\leq C_kt^{N_k}.
	\end{equation} By $s=\frac{t^3}{6}+O(t)$, this implies \begin{equation}\label{16}
		|\nabla_{g_L}^kw|_{g_L}\leq C_kt^{N_k}e^{-s}\leq C_ke^{-(1-\epsilon)s},\quad|\nabla_{g_L}^kw^{-1}|_{g_L}\leq C_kt^{N_k}e^s\leq C_ke^{(1+\epsilon)s}.
	\end{equation}

	Use \labelcref{15,16} to replace the covariant derivatives estimates of coordinates in \cite[Proposition 3.4]{heinNilpotentStructuresCollapsing2021}, the remaining steps are unchanged.
\end{proof}

\begin{remark}
	The holomorphic volume form on $M$ in this proposition is normalized by $\operatorname{Res}_D\Omega_M^{\mathrm{nor}}=\Omega_D$. For an arbitrary choice $\Omega_M=A\Omega_M^{\mathrm{nor}}$, we can start from the rescaled class $|A|^{-2/3}\mathfrak{c}$ and multiply the resulting Calabi--Yau metric by $|A|^{2/3}$. So from now on, we will always assume $\Omega_M=\Omega_M^{\mathrm{nor}}$.
\end{remark}

Take a defining section $S\in H^0(\overline{M},\mathcal{O}_{\overline{M}}(D))$, $\tau=-\log|S|_{\mathcal{O}_{\overline{M}}(D)}^2$, and hence $\alpha=\sqrt{-1}\partial\bar{\partial}\tau$. The property of exponential map gives \begin{equation}
	\Phi^*\tau=2s+O(|\xi|_L).
\end{equation} Together with \labelcref{st}, there is \begin{equation}\label{taut}
	t^3=3\Phi^*\tau+O(t).
\end{equation} In particular, $\tau$ is an exhaustion function on $M$.

Now we define the end form $\omega_{\mathrm{end}}=\beta+\sqrt{-1}\partial\bar{\partial}(U\circ\Phi^{-1})$ on a puncture neighborhood of $D$ in $\overline{M}$, where $U$ is the potential in \cref{potential}. Write $V=U\circ\Phi^{-1}$, then \begin{equation}\label{V}
	V=\frac{1}{4}(3\tau)^{4/3}+O(\tau).
\end{equation} The principal part is exactly the Tian--Yau metric's asymptotic in dimension $3$.

By \cref{cstr} we obtain \begin{equation}
	|\nabla_{g_L}^k(\Phi^*\omega_{\mathrm{end}}-\omega_L)|_{g_L}=O(e^{-ct^3}).
\end{equation} In particular, this implies $\omega_{\mathrm{end}}$ is K\"{a}hler in a sufficiently small puncture neighborhood of $D$ in $\overline{M}$. Also, \begin{equation}
	\log(\frac{\omega_{\mathrm{end}^3}}{6\sqrt{-1}\Omega_M\wedge\overline{\Omega}_M})=\log(\frac{\omega_t^2}{2h_t\cdot2\Omega_D\wedge\overline{\Omega}_D})+O(e^{-ct^3})=O(t^{-5}).
\end{equation}

Fix a small $\epsilon>0$, choose a smooth, convex and even function $\eta_\epsilon$ on $\mathbb{R}$, $\eta_\epsilon(x)=|x|$ if $|x|>\epsilon$, and that $|\eta_\epsilon'|\leq1$. Define a smooth maximum function by \begin{equation}
	\widetilde{\max}_\epsilon(x,y)=\frac{1}{2}(x+y+\eta_\epsilon(x-y))=\left\{\begin{aligned}
		&x,&&x-y>\epsilon \\
		&y,&&x-y<-\epsilon.
	\end{aligned}\right.
\end{equation}

Now we glue the end metric $\omega_{\mathrm{end}}=\beta+\sqrt{-1}\partial\bar{\partial}V$ with the background K\"{a}hler metric $\kappa=A\alpha+\beta=\beta+\sqrt{-1}\partial\bar{\partial}A\tau$.

First we take $T_0>0$ such that $\omega_{\mathrm{end}}$ is K\"{a}hler on $\{\tau>T_0-1\}$, and find \begin{equation}\label{14}
	C_0>\max_{\tau=T_0}(V-A\tau)+2\epsilon.
\end{equation} Since \begin{equation}
	\min_{\tau=T}(V-A\tau)\to+\infty
\end{equation} as $T\to+\infty$, there exists $T_1>T_0$ such that \begin{equation}
	\min_{\tau=T_1}(V-A\tau)>C_0+2\epsilon.
\end{equation} Then we define \begin{equation}\label{globalglue}
	\omega_0=\left\{\begin{aligned}
		&\omega_{\mathrm{end}},&&\tau>T_1 \\
		&\beta+\sqrt{-1}\partial\bar{\partial}\widetilde{\max}_\epsilon(V,A\tau+C_0),&&T_0\leq\tau\leq T_1 \\
		&\kappa,&&\tau<T_0.
	\end{aligned}\right.
\end{equation} By openness there is a small neighborhood of $\tau=T_0$ such that $V<A\tau+C_0$, and a small neighborhood of $\tau=T_1$ such that $V>A\tau+C_0$. So $\omega_0$ is a smooth K\"{a}hler metric on $M$.

\subsection{The integral condition}

Note the undetermined additive constant $c$ in the relation between $s$ and $t$ in \labelcref{st}. Changing this constant rescales the fiber coordinates by a constant multiple. We can adjust $c$ to require \begin{equation}
	\int_M(\omega_0^3-6\sqrt{-1}\Omega_M\wedge\overline{\Omega}_M)=0.
\end{equation} This is exactly the way used in \cite{tianCompleteKahlerManifolds1990}.

We write $\omega_0(c)$ to indicate dependence of the construction in $c$. Since \begin{equation}
	\frac{\omega_0(c)^3}{6\sqrt{-1}\Omega_M\wedge\overline{\Omega}_M}=1+O(t^{-5}),
\end{equation} and the volume form $6\sqrt{-1}\Omega_M\wedge\overline{\Omega}_M$ is of order $O(t^2)dt$ in the $t$ direction, we know \begin{equation}
	I(c):=\int_M(\omega_0(c)^3-6\sqrt{-1}\Omega_M\wedge\overline{\Omega}_M)
\end{equation} is finite.

Although the construction \labelcref{globalglue} depends on the gluing threshold $T_0,T_1$, and there is not a uniform choice for every shifting $c\in\mathbb{R}$. However, when changing $T_0,T_1$, the resulting $\omega_0(c)$ differs from a compactly supported potential function. So $I(c)$ is well-defined and depends only on $c$.

\begin{proposition}
	There exists $c\in\mathbb{R}$ such that $I(c)=0$.
\end{proposition}

\begin{proof}
	Recall that $s=-\log\|w\|$ is fixed by our construction. The change of $c$ makes $t=t(c,s)$ solved from \cref{st}, which exists from inverse function theorem. By the definition of $U=U(c,s)$ in \cref{potential}, we compute \begin{equation}
		\frac{\partial}{\partial c}U(c,s)=-2t(c,s).
	\end{equation}

	Locally we can take a uniform gluing threshold, and understand $\frac{d}{dc}\omega(c)$ in this sense. The potential of $\frac{d}{dc}\omega(c)$ is $\frac{\partial V}{\partial c}$ on the end. Therefore, \begin{equation}
		\begin{aligned}
			\frac{d}{dc}I(c)&=\lim_{T\to+\infty}\int_{(\Phi^{-1})^*t\leq T} 3\frac{d}{dc}\omega_0(c)\wedge\omega_0(c)^2 \\
			&=\lim_{T\to+\infty}\int_{(\Phi^{-1})^*t=T}\frac{3}{2}d^c(\frac{\partial V}{\partial c})\wedge\omega_0(c)^2.
		\end{aligned}
	\end{equation}

	Since $\Phi^*J_{\overline{M}}$ and $J_L$ differ exponentially small, the error from transferring $M$ to $L^\times$ vanishes after taking $T\to+\infty$. So \begin{equation}
		\begin{aligned}
			\frac{d}{dc}I(c)&=-\lim_{T\to+\infty}\int_{t=T}3d^ct\wedge\omega_L(c)^2 \\
			&=-\lim_{T\to+\infty}\int_D6\pi h_T^{-1}\omega_T^2 \\
			&=-\int_D 12\pi\cdot2\Omega_D\wedge\overline{\Omega}_D \\
			&=-48\pi^3c_1(L)^2.
		\end{aligned}
	\end{equation} Integrate, then \begin{equation}
		I(c)=I(0)-48\pi^3c_1(L)^2 c
	\end{equation} and $c$ is solved.
\end{proof}

From now on, we denote $\omega_0=\omega_0(c)$ for the choice of $c$ such that $I(c)=0$.

\section{Existence of complete Ricci-flat metric}\label{s7}

It remains to pass from the global approximate metric to an exact
Ricci-flat metric. We first recall the version of the Tian--Yau--Hein
package used for this purpose, and then verify its hypotheses for
$\omega_0$. The estimates from the model construction must now be
interpreted in terms of the geometry of the complete end. Once this is
done, the improved error estimate and the integral normalization give a
solution of the complex Monge--Amp\`{e}re equation. We conclude by checking
that the resulting metric is complete and remains in the prescribed
K\"{a}hler class.

\subsection{Tian--Yau--Hein package}

We'll solve a complex Monge--Amp\`{e}re equation on $M$ to perturb $\omega_0$ to the genuine Ricci-flat metric. We briefly recall the Tian--Yau--Hein package \cite{tianCompleteKahlerManifolds1990,heinGravitationalInstantons2010} used for solving this equation on noncompact manifold here.

\begin{definition}
	Let $\beta>0$, a complete noncompact Riemannian manifold $(M^n,g)$ is said to be $\mathrm{SOB}(\beta)$, if there exists $x_0\in M$, $r(x)=\dist(x_0,x)$, $C\geq1$ such that: \begin{enumerate}
		\item For any $s>C$, any points $p,q\in M$ with $r(p)=r(q)=s$ are connected by a path within the annulus $A(x_0,\frac{1}{C}s,Cs)$;
		\item For any $s>C$, $|B(x_0,s)|\leq Cs^\beta$;
		\item For any $x$ with $r(x)>C$, $|B(x,(1-\frac{1}{C})r(x))|\geq\frac{1}{C}r(x)^\beta$;
		\item For any $x$ with $r(x)>C$, $\Ric(x)\geq-Cr(x)^{-2}$.
	\end{enumerate}
\end{definition}

\begin{definition}
	Let $0\leq\lambda\leq1$, $k\geq0$, $0<\alpha<1$, a complete K\"{a}hler manifold $(M^n,\omega,g)$ is said to be $\mathrm{HMG}(\lambda,k,\alpha)$, if there exists $x_0\in M$, $r(x)=\dist(x_0,x)$, $C\geq1$ such that: \begin{enumerate}
		\item For any point $x$ with $r(x)\geq C$, there exists a holomorphic local diffeomorphism $\Phi_x$ from the unit ball $B\subset\mathbb{C}^n$ into $M$ such that $\Phi_x(0)=x$ and $\Phi_x(B)\supset B(x,\frac{1}{C}r(x)^\lambda)$;
		\item $h=r(x)^{-2\lambda}\Phi_x^*g$ satisfies $\operatorname{injrad}(h)\geq\frac{1}{C}$, $\frac{1}{C}g_{\mathrm{Euc}}\leq h\leq Cg_{\mathrm{Euc}}$ and $\|h-g_{\mathrm{Euc}}\|_{C^{k,\alpha}(B,g_{\mathrm{Euc}})}\leq C$.
	\end{enumerate}
\end{definition}

\begin{lemma}\label{hmg}
	A complete K\"{a}hler manifold with $|\mathrm{Rm}|+\sum_{i=1}^kr^{\lambda i}|\nabla^i\mathrm{Scal}|\leq Cr^{-2\lambda}$ is $\mathrm{HMG}(\lambda,k+1,\alpha)$ for every $0<\alpha<1$.
\end{lemma}

An $\mathrm{HMG}(\lambda,k,\alpha)$ structure induces a global H\"{o}lder norm on $M$ by \begin{equation}
	\|u\|_{C^{l,\gamma}(M)}:=\|u\|_{C^{l,\gamma}(B(x_0,2C),g)}+\sup_{r(x)\geq C}\|u\circ\Phi_x\|_{C^{l,\gamma}(B,g_{\mathrm{Euc}})}.
\end{equation}

\begin{theorem}\label{tyh}
	Let $(M^n,\omega,g)$ be a complete noncompact K\"{a}hler manifold that is $\mathrm{SOB}(\beta)$ and $\mathrm{HMG}(0,3,\alpha)$ for some $\beta>0$, $0<\alpha<1$. Let $r$ be a distance function to a fixed point. Let $f\in C^{2,\alpha}(M)$ satisfy $|f|\leq Cr^{-2-\epsilon}$ on $\{r>1\}$ for some $\epsilon>0$. If $\beta\leq2$, then assume in addition $\int_M(e^f-1)\omega^n=0$. Then there exists $\bar{\alpha}\in(0,\alpha]$ and $u\in C^{4,\bar{\alpha}}(M)$ satisfying \begin{equation}
		(\omega+\sqrt{-1}\partial\bar{\partial}u)^n=e^f\omega^n.
	\end{equation} If $\beta\leq2$, then moreover $\int_M|\nabla u|^2\omega^n<+\infty$. Independent of the value of $\beta$, if in addition $f\in C_{\mathrm{loc}}^{k,\bar{\alpha}}(M)$ for some $k\geq3$, then all such solutions $u$ belong to $C_{\mathrm{loc}}^{k+2,\bar{\alpha}}(M)$.
\end{theorem}

\subsection{Riemannian geometry of $\omega_0$}

The next four propositions verify the geometric hypotheses of
\cref{tyh}. We first compare the moment map coordinate with geodesic
distance and then estimate the Ricci tensor. The volume and connectivity
estimates establish $\mathrm{SOB}(\frac{3}{2})$, while the curvature estimates give
the required $\mathrm{HMG}$ condition. The distance comparison also
converts the decay obtained in \cref{s5} into the decay rate required by
the existence theorem.

\begin{proposition}\label{r}
	Fix a base point, then the distance function \begin{equation}
		r\sim t^2
	\end{equation} on the end. In particular, $(M,\omega_0)$ is complete.
\end{proposition}

\begin{proof}
	On the model $L$, $|dt|_{g_L}=h^{-1/2}\leq Ct^{-1}$. So $|d(t^2)|_{g_L}\leq C$. Integrating along a curve from the base point we get \begin{equation}
		r\geq\frac{1}{C}\int_{x_0}^xd(t^2)\geq\frac{1}{C}t^2.
	\end{equation}

	For the upper bound, first move on the slice $\{t=t(x_0)\}$, and then move along the radial $t$-curve. The diameter of the slice $\{t=t(x_0)\}$ is fixed. And since $h\sim t^2$, the radial length is bounded by $\int tdt\leq Ct^2$.

	Since $\Phi^*\omega_{\mathrm{end}}$ differs from $\omega_L$ exponentially small, the same comparison holds for $\omega_0$.
\end{proof}

\begin{proposition}\label{ric}
	We have \begin{equation}
		|\Ric(\omega_0)|\leq Cr^{-5/2}
	\end{equation} on the end.
\end{proposition}

\begin{proof}
	It's sufficient to prove on the model. Let \begin{equation}
		\rho=\log(\frac{\omega_t^2}{2h_t\cdot2\Omega_D\wedge\overline{\Omega}_D}),\quad \Ric(\omega_L)=-\sqrt{-1}\partial\bar{\partial}\rho.
	\end{equation} From \cref{metric} we know $\|\rho\|_{C_{\delta-2,t}^{k,\alpha;m}}\leq Ct^{-5}$. Apply \labelcref{eddc} we know \begin{align}
		&|d_Dd_D^c\rho|_{\omega_t}\leq Ct^{-5}l_t^{\delta-4}\leq Ct^{-5}, \\
		&h_t^{-1/2}|d_D\partial_t\rho|_{\omega_t}+h_t^{-3/2}|\partial_t\rho||d_Dh_t|_{\omega_t}\leq Ct^{-7}l_t^{\delta-3}\leq Ct^{-7}, \\
		&h_t^{-1}|\partial_t\rho|\cdot|\partial_t\omega_t|_{\omega_t}+h_t^{-1}|\partial_t^2\rho|+h_t^{-2}|\partial_t\rho||\partial_th_t|\leq Ct^{-9}l_t^{\delta-2}\leq Ct^{-9}.
	\end{align} And therefore \begin{equation}
		|\Ric(\omega_L)|\leq Ct^{-5}\leq Cr^{-5/2}.
	\end{equation}
\end{proof}

\begin{proposition}
	$(M,\omega_0)$ is $\mathrm{SOB}(\frac{3}{2})$.
\end{proposition}

\begin{proof}
	On the model, the volume of the slab $\{T\leq t\leq2T\}\subset L^\times$, \begin{equation}
		|\{T\leq t\leq2T\}|=\int_{T\leq t\leq2T}\frac{1}{2}\mathrm{d}t\wedge\Theta\wedge\omega_t^2\leq C\int_T^{2T}t^2dt\leq CT^3.
	\end{equation} So the volume growth $|B(x_0,R)|\leq CR^{3/2}$.

	For the relative lower bound, we first show $B_L(x,\frac{1}{2}r(x))$ contains the slab $\{(1-\epsilon)t(x)\leq t\leq(1+\epsilon)t(x)\}$. Since the gluing metric $\hat{\omega}_a$ in \cref{3.2} has a uniform diameter, \begin{equation}
		\diam(D,\omega_T)\sim T^{1/2}.
	\end{equation} Lift the shortest geodesic to the circle bundle $(\{t=T\},g_T+h_T^{-1}\Theta^2)$, the circle fiber length is bounded by $2\pi h_T^{-1/2}\leq CT^{-1}$. So the diameter of the circle bundle remains \begin{equation}
		\diam(\{t=T\},g_T+h_T^{-1}\Theta^2)\sim T^{1/2}.
	\end{equation} Now we take $y$ belonging to the slab $(1-\epsilon)t(x)\leq t(y)\leq(1+\epsilon)t(x)$. The radial distance between $x,y$ is controlled by $2C\epsilon t(x)^2$, and the slice distance by $Ct(x)^{1/2}$, hence \begin{equation}
		\dist(x,y)\leq C(\epsilon t(x)^2+t(x)^{1/2})\leq\frac{1}{2}r(x)
	\end{equation} once $\epsilon$ is small enough and $x$ lies sufficiently far out on the end.

	Since $\nabla t=h_t^{-1}\partial_t$, by coarea formula, \begin{equation}
		\begin{aligned}
			|B_L(x,\frac{1}{2}r(x))|&\geq\int1_{(1-\epsilon)t(x)\leq t\leq(1+\epsilon)t(x)}\frac{\omega_L^3}{6} \\
			&\geq \frac{t(x)}{C}\int1_{(1-\epsilon)t(x)\leq t\leq(1+\epsilon)t(x)}|\nabla t|\frac{\omega_L^3}{6} \\
			&=\frac{t(x)}{C}\int_{(1-\epsilon)t(x)}^{(1+\epsilon)t(x)}\operatorname{Area}(\{t=\tau\})d\tau \\
			&\geq \frac{1}{C}t(x)^3,
		\end{aligned}
	\end{equation} where the area of the level set \begin{equation}
		\operatorname{Area}(\{t=\tau\})=\int_D 2\pi h_\tau^{-1/2}\frac{\omega_\tau^2}{2}\sim\tau.
	\end{equation} Transfer to $M$, then $|B(x,\frac{1}{2}r(x))|\geq\frac{1}{C}r(x)^{3/2}$.

	The Ricci decay is proved in \cref{ric}. Finally, consider the connectedness condition. Let $x,y$ with $r(x)=r(y)=R$ on the model. Then their $t$-coordinates differ from $CR^{1/2}$. We use radial segment to transport them in the same $t$-coordinate, and connect by slice geodesic controlled by $\diam(\{t=T\})\leq CR^{1/4}$. Clearly, this path is contained in some $A(x_0,\frac{1}{C}R,CR)$.
\end{proof}

\begin{proposition}
	$(M,\omega_0)$ is $\mathrm{HMG}(0,k,\alpha)$ for $k\geq3$, $0<\alpha<1$.
\end{proposition}

\begin{proof}
	We use \cref{hmg}. The Riemannian metric associated to $\omega_L$ is \begin{equation}
		g_L=h_t^{-1}\Theta^2+h_tdt^2+g_t.
	\end{equation} View $L$ as a circle bundle over $B=D\times\mathbb{R}_{\geq0}$, and let \begin{equation}
		g_B=h_tdt^2+g_t
	\end{equation} be the base metric.

	On $B$, $\nabla t=h_t^{-1}\partial_t$, $|\nabla t|^2=h_t^{-1}$. So the unit vector field $\nu_t=h_t^{-1/2}\partial_t$. We compute the second fundamental form on the slice \begin{equation}
		\mathrm{II}=h_t^{1/2}\Hess t=h_t^{1/2}\cdot\frac{1}{2}\mathcal{L}_{\nabla t}g_B=\frac{1}{2}h_t^{-1/2}\partial_tg_t.
	\end{equation} Therefore, \begin{equation}
		|\mathrm{II}|_{g_B}+|\nabla\mathrm{II}|_{g_B}+|\partial_t\mathrm{II}|_{g_B}\leq Ct^{-2}
	\end{equation} from the properties of $(\omega_t,h_t)$. Denote $i,j,k,l$ the tangential directions along the slice, $t$ the radial direction $\partial_t$, then \begin{align}
		&(\mathrm{Rm}_B)_{ijkl}=(\mathrm{Rm}_{g_t})_{ijkl}-\mathrm{II}_{il}\mathrm{II}_{jk}+\mathrm{II}_{ik}\mathrm{II}_{jl}, \\
		&(\mathrm{Rm}_B)_{ijkt}=-h_t^{1/2}\nabla_i\mathrm{II}_{jk}+h_t^{1/2}\nabla_j\mathrm{II}_{ik}, \\
		&(\mathrm{Rm}_B)_{itjt}=h_t^{1/2}\partial_t\mathrm{II}_{ij}-h_tg_t^{kl}\mathrm{II}_{ik}\mathrm{II}_{jl}+h_t^{1/2}\nabla_i\nabla_jh_t^{1/2}.
	\end{align} The first two are the Gauss--Codazzi equations, the last is from the direct computation: \begin{equation}
		\begin{aligned}
			&\Gamma_{tt}^t=\frac{1}{2}h_t^{-1}\partial_th_t,\quad\Gamma_{tj}^t=\frac{1}{2}h_t^{-1}\partial_jh_t, \\
			&\Gamma_{ij}^t=-\frac{1}{2}h_t^{-1}\partial_tg_{ij}=-h_t^{-1/2}\mathrm{II}_{ij},\quad\Gamma_{tj}^k=\frac{1}{2}g^{kl}\partial_t g_{jl}=h_t^{1/2}g^{kl}\mathrm{II}_{jl},
		\end{aligned}
	\end{equation} and \begin{equation}
		(\mathrm{Rm}_B)_{itjt}=h_t(\partial_i\Gamma_{tj}^t-\partial_t\Gamma_{ij}^t+(\Gamma_{tj}^k\Gamma_{ik}^t+\Gamma_{tj}^t\Gamma_{it}^t-\Gamma_{ij}^k\Gamma_{tk}^t-\Gamma_{ij}^t\Gamma_{tt}^t)).
	\end{equation} From \cref{qatlas}, and the coefficients condition of $\omega_t$, we obtain \begin{equation}
		|\mathrm{Rm}_{g_t^{(0)}}|\leq Cl_t^{-2}\leq C,\quad|\mathrm{Rm}_{g_t}|\leq C.
	\end{equation} Therefore, \begin{equation}
		|\mathrm{Rm}_B|\leq C.
	\end{equation}

	Now the curvature of the circle bundle's connection \begin{equation}
		F_\Theta=d\Theta=d_D\Theta+dt\wedge\partial_t\Theta=\partial_t\omega_t-dt\wedge d_D^ch_t,
	\end{equation} so that \begin{equation}
		|F_\Theta|+|\nabla^BF_\Theta|\leq Ct^{-1}.
	\end{equation} The calculation in \cite[Lemma 3.2]{shankarNonnegativelyPositivelyCurved2005} implies \begin{equation}
		|\mathrm{Rm}_L|\leq C(|\mathrm{Rm}_B|+h_t^{-1}|F_\Theta|^2+h_t^{-1/2}|\nabla^BF_\Theta|+|d(h_t^{-1/2})||F_\Theta|+h^{1/2}|\Hess_B(h^{-1/2})|)\leq C.
	\end{equation}

	Meanwhile, \begin{equation}
		\mathrm{Scal}_L=-\Delta_L\rho,\quad\rho=\log(\frac{\omega_t^2}{2h_t\cdot2\Omega_D\wedge\overline{\Omega}_D}).
	\end{equation} For an $S^1$-invariant tensor $T$, compute in local coordinates, we can estimate \begin{equation}\label{tensor}
		|\nabla_{g_L}T|_{g_L}\leq C(|\nabla_{g_t}\mathbf{T}|_{g_L}+t^{-1}|\partial_t\mathbf{T}|_{g_L}+t^{-2}|T|_{g_L}),
	\end{equation} where $\mathbf{T}$ is the collection of $T$'s components according to the decomposition \begin{equation}
		TL=(TD)^H\oplus\mathbb{R}(h_t^{-1/2}\partial_t)\oplus\mathbb{R}(h_t^{1/2}\partial_\sigma)
	\end{equation} and the differentials apply to every component. Since $\|\rho\|_{C_{\delta-2,t}^{k,\alpha;m}}\leq Ct^{-5}$ from \cref{metric}, there is \begin{equation}
		|\nabla_{g_L}^{i+2}\rho|\leq Ct^{-5},\quad|\nabla_{g_L}^i\mathrm{Scal}_L|\leq Ct^{-5}.
	\end{equation} In summary, we have \begin{equation}
		|\mathrm{Rm}_L|+\sum_{i=1}^k|\nabla_{g_L}^i\mathrm{Scal}_L|\leq C.
	\end{equation} This shows $\omega_L$ is $\mathrm{HMG}(0,k+1,\alpha)$, so is $\omega_0$.
\end{proof}

Apply the Tian--Yau--Hein package to solve the equation \begin{equation}
	(\omega_0+\sqrt{-1}\partial\bar{\partial}u)^3=e^f\omega_0^3,\quad f=\log(\frac{6\sqrt{-1}\Omega_M\wedge\overline{\Omega}_M}{\omega_0^3}).
\end{equation} \Cref{metric} and \labelcref{tensor} imply $f\in C^{2,\alpha}(M)$ and $|f|\leq Cr^{-5/2}$. All hypotheses of \cref{tyh} are therefore satisfied. Applying the theorem yields a smooth function $u$ such that \begin{equation}
	\omega_{\mathrm{CY}}=\omega_0+\sqrt{-1}\partial\bar{\partial}u>0,\quad\omega_{\mathrm{CY}}^3=6\sqrt{-1}\Omega_M\wedge\overline{\Omega}_M.
\end{equation} The global $C^{4,\bar{\alpha}}(M)$ bound for $u$ gives a uniform upper bound for the eigenvalues of $\omega_{\mathrm{CY}}$ relative to $\omega_0$. Since their product is $e^f$, which is bounded away from $0$, the eigenvalues also have a uniform positive lower bound. Thus $\omega_{\mathrm{CY}}$ is uniformly equivalent to $\omega_0$ and is complete. Finally, \begin{equation}
	[\omega_{\mathrm{CY}}]=[\omega_0]=[\beta|_M]=\mathfrak{c},
\end{equation} so the metric lies in the prescribed K\"{a}hler class.

\begin{theorem}
	Let $\overline{M}$ be a smooth compact K\"{a}hler threefold, $D\in |-K_{\overline{M}}|$ be a smooth K3 surface, $M=\overline{M}\setminus D$, $\Omega_{M}$ be a meromorphic volume form with $\operatorname{div}\Omega_{M}=-D$. Assume the normal bundle of $D$ is nef and big. Then in every K\"{a}hler class of $M$ such that all the exceptional periods are positive, there exists a complete Ricci-flat K\"{a}hler metric $\omega_{\mathrm{CY}}$ on $M$, such that $\omega_{\mathrm{CY}}^3=6\sqrt{-1}\Omega_{M}\wedge\overline{\Omega}_{M}$.
\end{theorem}

\bibliographystyle{ytamsalpha}
\bibliography{ref}

\end{document}